\documentclass[preprint,12pt]{elsarticle}
\usepackage{listings}
\usepackage{placeins}
\usepackage{amsmath}
\usepackage[numbers]{natbib}
\usepackage{algorithm,algorithmicx}
\usepackage{physics}

\usepackage[english]{babel} 

\usepackage{epsfig}

\usepackage{algorithmicx}
\usepackage{algpseudocode}

\usepackage{stmaryrd}
\usepackage{amsmath} 
\usepackage{color}
\usepackage{amssymb} 
\usepackage{amsfonts} 
\usepackage{mathrsfs}

\usepackage[T1]{fontenc}
\usepackage[utf8]{inputenc}
\usepackage[font=small,labelfont=bf,tableposition=top]{caption}
\usepackage{multicol}

\usepackage{subcaption}
\usepackage{graphicx}
\input{epsf}

\newcommand{\comment}[1] {}

\newcommand{\diff}[2] {\frac{\partial #1}{\partial #2}}

\newcommand{\difft}[2] {\frac{d #1}{d #2}}

\newcommand{\vc}[1]{\mbox{\boldmath$#1$\unboldmath}}
\newcommand{\inte}{\int_{\mbox{\footnotesize ${\Omega_e}$}}}

\newcommand{\intce}{\int_{\mbox{\footnotesize ${\widehat{\Omega}}$}}}

\newcommand{\wh}{\widehat}

\newcommand{\be}{\begin{equation}}
\newcommand{\ee}{\end{equation}}
\newcommand{\bea}{\begin{eqnarray*}}
\newcommand{\eea}{\end{eqnarray*}}

\newcommand{\transpose}{^{\mathcal{T}}}

\newcommand{\qvector}{\vc{q}}

\newcommand\ST{\rule[-0.75em]{0pt}{2em}}

\newcommand{\Lfunction}{\mathbb{L}}

\newcommand{\Efunction}{E}
\newcommand{\Ifunction}{I}

\newcommand{\gmo}{\left(\gamma-1\right)} 

\newcommand{\ignore}[1]{}

\begin{document}
\begin{frontmatter}

\title{A Performance Study of Horizontally Explicit Vertically Implicit (HEVI) Time-Integrators for Non-Hydrostatic Atmospheric Models}

\author[1]{Francis X. Giraldo\corref{cor1}}
\ead{fxgirald@nps.edu}

\author[1]{Felipe A.\ V.\ de Bragan\c{c}a Alves\corref{cor2}}
\ead{favbalves@gmail.com}

\author[2]{James F. Kelly}
\ead{james.kelly@nrl.navy.mil}

\author[1]{Soonpil Kang}
\ead{soonpil.kang.ks@nps.edu}

\author[3]{P. Alex Reinecke}
\ead{alex.reinecke@nrlmry.navy.mil}

\cortext[cor1]{Corresponding author}
\cortext[cor2]{Now at Instituto de Matem\'{a}tica e Estat\'{i}stica, Universidade de S\~{a}o Paulo, Brazil}

\affiliation[1]{organization={Department of Applied Mathematics, Naval Postgraduate School},
            city={Monterey},
            state={CA},
            country={U.S.}}
\affiliation[2]{organization={Space Science Division, U.S. Naval Research Laboratory},
            city={Washington, DC},
            country={U.S.}}
\affiliation[3]{organization={Marine Meteorology Division, U.S. Naval Research Laboratory},
            city={Monterey, CA},
            country={U.S.}}

\begin{abstract}
We conduct a thorough study of different forms of horizontally explicit and vertically implicit (HEVI) time-integration strategies for the compressible Euler equations on spherical domains typical of nonhydrostatic global atmospheric applications. We compare the computational time and complexity of two nonlinear variants (NHEVI-GMRES and NHEVI-LU) and a linear variant (LHEVI). We report on the performance of these three variants for a number of additive Runge-Kutta Methods ranging in order of accuracy from second through fifth, and confirm the expected order of accuracy of  the HEVI methods for each time-integrator. To gauge the maximum usable time-step of each HEVI method, we run simulations of a nonhydrostatic baroclinic instability for 100 days and then use this time-step to compare the time-to-solution of each method. The results show that NHEVI-LU is 2x faster than NHEVI-GMRES, and LHEVI is 5x faster than NHEVI-LU, for the idealized cases tested. 
The baroclinic instability and inertia-gravity wave simulations indicate that the optimal choice of time-integrator is LHEVI with either second or third order schemes, as both schemes yield similar time to solution and relative L2 error at their maximum usable time-steps.
In the future, we will report on whether these results hold for more complex problems using, e.g., real atmospheric data and/or a higher model top typical of space weather applications.
\end{abstract}


\end{frontmatter}

\section{Introduction}
\label{sec:introduction}

All problems in the numerical solution of time-dependent partial differential equations (PDEs) require, at the very least, two key ingredients for accurate, stable, and efficient solutions. These ingredients are the selection of the methods for approximating the (i) spatial and (ii) temporal (time) derivatives.  The focus of this paper is on time-integration methods specifically designed for geophysical fluid dynamics (GFD) where the stiffness of the PDEs results from a finer grid spacing in the vertical relative to the horizontal dimension.  We are interested in solving the nonhydrostatic form of the compressible Euler equations, where explicit time-integration is impractical due to the small time-steps dictated by the acoustic speed and reflected by the Courant-Friedrichs-Lewy condition \cite{courant:1967}. For this reason, recent work has focused on nonlinear horizontally explicit and vertically implicit methods (HEVI) where, as the name suggests, the horizontal flow along the surface of the earth is handled explicitly in time, whereas the flow normal to the surface of the earth along the radial (or vertical) direction is handled implicitly. The work by Gardner et al.\ \cite{gardner:2018} describes a nonlinear HEVI method that uses a GMRES-based Krylov solver, which we call nonlinear HEVI (or NHEVI-GMRES). Since the publication of that paper, others have contributed to this discussion and we mention in particular the work by Steyer et al.\ \cite{steyer:2019} and Vogl et al.\ \cite{vogl:2019}. The goal of our current work is to analyze NHEVI and seek improvements to it such that (i) the stability and accuracy of the original method is not degraded and (ii) its efficiency is enhanced. Some researchers have reached the conclusion that NHEVI-GMRES may be too expensive and modifications have been studied but perhaps not published (e.g., work on the U.S. Navy's NEPTUNE model), while other approaches have already appeared in the literature such as LHEVI (e.g., see \cite{giraldo:2013,weller:2013,lock:2014,bao:2015,baldauf:2021,lee:2021,sridhar:2022,waruszewski:2022,souza:2023} for IMEX-LHEVI and \cite{satoh:2002,klemp:2018} for split-explicit LHEVI) but in these works it was not compared to the NHEVI method.  Therefore, the motivation for the current work is to quantitatively compare the NHEVI method relative to linear HEVI (LHEVI); where we use the Nonhydrostatic Unified Model of the Atmosphere (NUMA) \cite{giraldo:2013} for the numerical experiments.

NUMA is a research CFD code that uses Element-Based Galerkin (EBG) methods \cite{marras:2015c,giraldo:2020} for spatial discretization and a suite of explicit, IMEX, and fully implicit time-integrators.  NEPTUNE is primarily a low-altitude NWP model, although we are also developing a high-altitude, whole atmosphere version \cite{akmaev:2011} of both NUMA and NEPTUNE for space weather applications.  Early versions of NUMA \cite{giraldo:2013} and  NEPTUNE \cite{reinecke:2016} relied primarily on 1D linear IMEX schemes, where each column of fluid was rotated to be aligned with the north pole during the implicit step.  In addition, linear IMEX relies on an \emph{a priori} reference state to construct the linear operator.  Both of these constraints are problematic for an NWP model.  The rotations add additional computational cost to the solver, while pre-defined reference states are typically not available for simulations with real data.  In addition, in whole atmosphere applications, the upper atmosphere undergoes very large day/night transitions in temperature (hundreds of degrees Kelvin), which makes a predefined reference state numerically unstable.  For these reasons, we decided to implement and 
evaluate various HEVI schemes for both the low-altitude and high-altitude versions of NUMA.  The complexity and numerical efficiency of each of these schemes are evaluated along with their accuracy with idealized benchmarks, allowing us to determine the optimal configuration for NUMA and then make a recommendation for NEPTUNE.


The remainder of the paper is organized as follows.
In Sec.\ \ref{sec:governing_equations}, we describe the form of the governing equations used in our study along with the contravariant formulation used in its discretization. Section \ref{sec:space} details the spatial discretization scheme based on EBG methods \cite{giraldo:2020}, including special treatment of the metric terms using curl-invariant metrics \cite{kopriva:2006}.  Section \ref{sec:hevi}  summarizes the three different flavors of the HEVI methods: nonlinear HEVI based on an iterative Krylov solver (NHEVI-GMRES), nonlinear HEVI using a direct LU solver (NHEVI-LU), and linear HEVI (LHEVI).  Section \ref{sec:complexity} provides a complexity analysis of the three HEVI variants, while Sec.\ \ref{sec:time_integrators}  compares five additive Runge-Kutta (ARK) time-integrators that may be used with any of the three HEVI variants.  Numerical results using two standard global atmospheric test cases are shown in Sec.\ \ref{sec:results}, followed by conclusions in Sec.\ \ref{sec:conclusions}.  In  \ref{sec:appendix/set2c_set3c} we discuss two other forms of the nonhydrostatic equations found in the literature and include the necessary ingredients for solving them using the HEVI time-integration methods described in this work.



\section{Governing Equations}
\label{sec:governing_equations}
In this study we use the compressible Euler equations (with an artificial hyper-diffusion operator). According to the convention we defined in \cite{giraldo:2008a,giraldo:2010b} we refer to this particular form of the Euler equations as set 2NC, where the "NC" refers to a non-conservative set. In \ref{sec:appendix/set2c_set3c/governing_equations}, we include two other equation sets: 2C which is the conservation form of 2NC and 3C which uses the total energy as a prognostic variable.
All three sets are valid for a deep-atmosphere \cite{white:2005}. In a forthcoming paper, we present an equation set which uses internal energy for the thermodynamic variable \cite{kelly:2023} that is valid for space weather applications.

\subsection{Set 2NC}
This is perhaps the most widely used form of the equations in atmospheric sciences and is written as follows
\begin{subequations}
\label{eq:set2nc}
\begin{equation}
\label{eq:set2nc_mass}
\diff{\rho}{t} + \nabla \cdot \left(\rho \vc{u} \right)= 0
\end{equation}
\begin{equation}
\label{eq:set2nc_momentum}
\diff{\vc{u}}{t} + \vc{u} \cdot \nabla  \vc{u} + \frac{1}{\rho} \nabla P + \nabla \Phi +  2 \omega \hat{\vc{r}} \times \vc{u} = \vc{0}
\end{equation}
\begin{equation}
\label{eq:set2nc_energy}
\diff{\theta}{t} + \vc{u} \cdot \nabla  \theta = 0
\end{equation}
\begin{equation}
\label{eq:set2nc_pressure}
P = P_A \left( \frac{ \rho R \theta}{P_A} \right)^{\gamma}.
\end{equation}
\end{subequations}
In these equations, $\rho$ is density, $\vc{u}=(u,v,w)^\mathcal{T}$ is the Cartesian velocity field, $\theta$ is the potential temperature, $P$ is pressure, $\Phi=g |\vc{r}| $ is geopotential, $\vc{r}$ is the position vector (measured from sea level), $\omega$ is the angular rotation of the earth, $\hat{\vc{r}}=\frac{\vc{r}}{|\vc{r}|}$ is the unit vector along the direction of $\vc{r}$, $R$ is the specific gas constant, $\gamma=\frac{c_p}{c_v}$ is the ratio of specific heats with respect to constant pressure and volume, $P_A=10^5$ Pa is the reference pressure at sea level, and the superscript $\mathcal{T}$ is the transpose operator. Note that the momentum \eqref{eq:set2nc_momentum} and thermodynamic \eqref{eq:set2nc_energy} equations are in non-conservation (advective) form while the continuity \eqref{eq:set2nc_mass} equation is in conservation (flux) form; with the proper numerics (as the EBG method that we use in this work) this equation set will globally conserve mass as long as we do not use the product rule in \eqref{eq:set2nc_mass}. Set 2NC, or variants thereof, are used in, e.g., E3SM \cite{taylor:2020} and NEPTUNE \cite{reinecke:2016}.

\subsection{Contravariant Formulations}
\label{sec:contravariant_form}
NUMA uses a Cartesian coordinate system, which while natural for flow in a box (where the grid is aligned with the coordinate axes), appears out of place for global atmospheric simulations on the sphere where the coordinate axes are no longer aligned with a spherical grid.   This issue can be ameliorated by using contravariant vector quantities that \emph{are} aligned with the grid.  This requires mapping the Cartesian coordinates to the reference coordinates such that vectors (e.g., velocity) are transformed as follows:
$(u,v,w) \rightarrow (u^{\xi},u^{\eta},u^{\zeta})$ where $(\xi,\eta)$ are aligned with the spherical manifold and $\zeta$ is aligned with the radial component, i.e., the direction along which the HEVI method is applied.

To express the transformation from covariant to contravariant vector quantities, we make use of Einstein notation whereby $u_i$ and $u^i$ denote covariant and contravariant quantities, respectively. In what follows, we refer to the covariant vectors as the vectors in Cartesian space which, within an element, represent a curvilinear representation, while the contravariant vectors are those that are aligned with the manifold and the crux of the problem is to construct a local (element-wise) map that takes each element defined in Cartesian space to its corresponding reference space.

It is easy to show (via the chain rule) that the advective terms in covariant form can be replaced by their contravariant form as follows
\be
\label{eq:grad}
\sum_{i=1}^3 u_i \diff{q}{x_i} = \sum_{i=1}^3 u^i \diff{q}{\xi^i} 
\ee
where $\vc{u}=(u,v,w)=(u_1,u_2,u_3)$, 
\be
\label{eq:grad_x}
\nabla=\left( \diff{}{x}, \diff{}{y}, \diff{}{z} \right)=\left( \diff{}{x_1}, \diff{}{x_2}, \diff{}{x_3} \right)
\ee
are the covariant vector quantities and 
$\vc{u}^{{\xi}}=\left( u^{\xi},u^{\eta},u^{\zeta} \right)=\left( u^1,u^2,u^3 \right)$, 
\be
\label{eq:grad_ksi}
\nabla_{{\xi}}=\left( \diff{}{\xi}, \diff{}{\eta}, \diff{}{\zeta} \right)=\left( \diff{}{\xi^1}, \diff{}{\xi^2}, \diff{}{\xi^3} \right)
\ee
 are the contravariant vector quantities. It is less easy to show that divergence (e.g., see \cite{kopriva:2009}) can be transformed as follows
\be
\label{eq:div_contra}
\sum_{i=1}^3 \diff{ }{x_i} \left(q u_i \right) = \frac{1}{J} \sum_{i=1}^3 \diff{ } {\xi^i} \left(J q u^i\right)
\ee
where $u^i= \vc{u} \cdot \nabla \xi^i$ for any vector $\vc{u}$ and scalar $q$ with metric Jacobian $J$ (see Eq.\ \eqref{eq:jacobian-metrics}) and contravariant metric terms $\nabla \xi^i$ given by \eqref{eq:curl-invariant-metrics} where $\xi^i$ is defined as $( \xi^1, \xi^2, \xi^3) = ( \xi, \eta, \zeta)$. With all these transformations defined, we can proceed to the spatial discretization of the equation sets.  
Note that discrete mass conservation requires the continuity equation \eqref{eq:set2nc_mass} to be discretized in the divergence form shown in \eqref{eq:div_contra}. We cannot use the product rule and expand the divergence operator but must maintain it in divergence form; it can be shown that the global integral of this term vanishes when using EBG methods (e.g., see \cite{taylor:2010,giraldo:2020})  that are presented in Sec.\ \ref{sec:space}. 
The same argument also holds for sets 2C and 3C described in \ref{sec:appendix/set2c_set3c/governing_equations} that are able to conserve quantities written in this form (e.g., the mass for both and density potential temperature and total energy, respectively).

\section{Spatial Discretization}
\label{sec:space}
Although the focus of this work is on time-integration, let us first describe the spectral element method used for the spatial discretization which will assist in describing the construction of the operator Jacobian in the HEVI methods.  The basis functions $\psi(\xi,\eta,\zeta)$ are constructed as a tensor-product $\otimes$ of the one-dimensional (1D) Lagrange polynomials $h(\xi^i)$ as follows
\be
\label{eq:space/basis}
\psi_{I}(\xi,\eta,\zeta)=h_i(\xi) \otimes h_j(\eta) \otimes h_k(\zeta)
\ee
where $I=i + j(N_{\xi}+1) + k(N_{\xi}+1)(N_{\eta}+1)$ is the map from the tensor-product to monolithic space,
with corresponding quadrature weights $\varpi_{I}=\varpi^{\xi}_i  \varpi^{\eta}_j  \varpi^{\zeta}_k$ with $(\xi_i,\eta_j,\zeta_k)$, $i=0,\ldots,N_{\xi}$, $j=0,\ldots,N_{\eta}$, $k=0,\ldots,N_{\zeta}$ denoting the Lobatto points of order $(N_{\xi},N_{\eta},N_{\zeta})$, respectively, that are used for both interpolation and integration, and $\psi \in \mathbb{H}^1$ (Sobolev space) (see, e.g., \cite{patera:1984,deville:2002,kopriva:2009,giraldo:2020}). Because we use the Lobatto points for both interpolation and integration, the 1D basis functions satisfy  $h_i(\xi^k_j) =\delta_{ij}$ for any of the contravariant coordinates $\xi^k$ where $\delta$ is the Kronecker function; this is due to the \emph{Cardinality property} of Lagrange polynomials.  Differentiating \eqref{eq:space/basis} yields
\begin{align}
\label{eq:space/basis/derivatives}
\diff{ }{\xi} \psi_{I}(\xi,\eta,\zeta) &= \difft{h_i(\xi)}{\xi} \otimes h_j(\eta) \otimes h_k(\zeta) \notag \\
\diff{ }{\eta} \psi_{I}(\xi,\eta,\zeta) &= h_i(\xi) \otimes \difft{h_j(\eta)}{\eta} \otimes h_k(\zeta)  \\
\diff{ }{\zeta} \psi_{I}(\xi,\eta,\zeta) &= h_i(\xi) \otimes h_j(\eta) \otimes \difft{h_k(\zeta)}{\zeta} \notag
\end{align}
where, e.g., we define the 1D differentiation matrices as $D_{ij}^{\xi^k} = \difft{h_i(\xi^k_j)}{\xi^k}$, $\forall \, i,j=0,\ldots,N_{\xi^k}$. To simplify the exposition in Sec.\ \ref{sec:hevi}, we exploit the \emph{row-sum property} of $D^{\xi^k}$ defined as $\sum_{i=0}^{N_{\xi^k}} D_{ij}^{\xi^k} = 0$ that can be used to show that for two functions $p$ and $q$ we satisfy the identity
\be
\label{eq:space/flux_difference}
p_i D_{ij}^{\xi^k} q_j = D_{ij}^{\xi^k} \left( p_i q_j \right),
\ee
which means that we can write advective terms, such as $u \diff{u}{\zeta}$ as $D^{\zeta}_{ij} \left( u_i u_j \right)$.
This \emph{flux-differencing} approach is used extensively in the literature addressing the construction of kinetic-energy preserving and entropy-stable discontinuous Galerkin methods (e.g., see \cite{gassner:2016,waruszewski:2022,souza:2023}). 

Using the basis functions \eqref{eq:space/basis} we can construct the mass matrix on each element $\Omega_e$ as such
\be
\label{eq:space/mass_matrix_element}
M^{(e)}_{ij} = \inte \psi_{i}(x,y,z) \psi_{j}(x,y,z) d\Omega_e \equiv \intce \psi_{i}(\xi,\eta,\zeta) \psi_{j}(\xi,\eta,\zeta) J d\wh{\Omega}
\ee
where we now use the monolithic space for the subscripts ($i,j$), $J$ is the metric Jacobian defined in the next section and $\wh{\Omega}$ is the reference element defined by the cube $(\xi,\eta,\zeta) \in [-1,+1]^3$; note that we never need the basis functions $\psi(x,y,z)$ defined in the physical coordinates, only those in the reference coordinates.
Applying co-located Lobatto integration to \eqref{eq:space/mass_matrix_element} yields the diagonal matrix
\be
\label{eq:space/mass_matrix_element/integration}
M^{(e)}_{ij} =  \varpi_i J_i \delta_{ij}
\ee
due to Cardinality but is inexact since $N$ Lobatto points integrate $\order{2N-1}$ polynomials exactly (the mass matrix represents a $2N$ polynomial). To denote the diagonal nature of the mass matrix, we write it as $M^{(e)}_{i}$.  To construct the global mass matrix we apply \emph{direct-stiffness summation} (DSS) as follows
\be
\label{eq:space/mass_matrix_element/integration}
M_{I} = \bigwedge_{e=1}^{N_e} M^{(e)}_i
\ee
where $(i,e) \rightarrow I$ denotes the map from local element-wise to global gridpoint notation (see \cite{kelly:2012,giraldo:2020}) and $N_e$ denotes the number of elements such that 
\[
\Omega = \bigcup_{e=1}^{N_e} \Omega_e
\]
where $\Omega$ is the global domain.
Below we drop the subscript and write the global gridpoint mass matrix as $M$ where its inverse is simply defined and denoted as $M^{-1}$.

\subsection{Metric Terms}
\label{sec:space/metric_terms}
To transform from the physical Cartesian coordinates $(x,y,z)$  to the reference element coordinates $(\xi,\eta,\zeta)$ requires us to discuss the construction of the metric terms. The standard approach to constructing them are the \emph{cross-product} metrics (e.g., see \cite{giraldo:2020}) written as 
\begin{equation}
\nabla \xi^i = \frac{1}{J} \left( \diff{\vc{x}}{\xi^j} \times \diff{\vc{x}}{\xi^k} \right)
\label{eq:cross-product-metrics}
\end{equation}
where the metric Jacobian is defined as 
\be
\label{eq:jacobian-metrics}
J=\diff{\vc{x} }{\xi} \cdot \left( \diff{\vc{x} }{\eta} \times \diff{\vc{x} }{\zeta} \right)
\ee
and $i,j,k$ are defined cyclically such that if $i=1$, then $j=2$, and $k=3$, etc.
The cross-product metric terms 
can maintain constant stream preservation in 2D but do not have this property in 3D (see \cite{kopriva:2006}).  To satisfy this condition in 3D requires the use of the \emph{curl-invariant} metrics (e.g., see \cite{kopriva:2006}) that are written as  
\begin{equation}
\nabla \xi^i = \frac{1}{2J} \left[   \diff{ }{\xi^k} \left( \diff{\vc{x}}{\xi^j} \times \vc{x} \right) - \diff{ }{\xi^j} \left( \diff{\vc{x}}{\xi^k} \times \vc{x} \right)  \right]
\label{eq:curl-invariant-metrics}
\end{equation}
where, once again, the indices $(i,j,k)$ are cyclic, and the metric Jacobian is computed as described in \cite{kopriva:2006}.

In previous work \cite{waruszewski:2022,souza:2023}, we used the curl-invariant metrics but did not compare them to the cross-product metrics. We rectify this by demonstrating the superiority of the curl-invariant metrics to the cross-product form and then use the curl-invariant metrics throughout the paper because they are guaranteed to maintain constant-stream-preservation under all conditions, including when the boundaries are under-resolved, whether topography is present or not, and under inexact integration of the integrals required in efficient EBG methods. In an upcoming paper, we quantitatively compare various metric terms using a suite of tropospheric and high-altitude atmospheric problems in \cite{kelly:2023}, further demonstrating the advantages of the curl-invariant metrics.

\subsection{Metric Terms along Columns}
\label{sec:space/metric_terms_columns}
\begin{figure}
\begin{center}
 \includegraphics[width=0.4\textwidth]{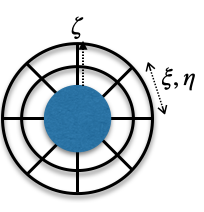}
 \caption{Alignment of grid such that $\xi-\eta$ define the spherical manifold while $\zeta$ defines the radial direction.}
 \label{fig:space/reference_grid}
\end{center}
\end{figure}

By constructing our grids in NUMA such that the $\zeta$ direction is aligned with the vertical columns, this means that the HEVI methods only have to consider the $\zeta$ direction; this is depicted in Fig.\ \ref{fig:space/reference_grid} where $\xi-\eta$ define the spherical manifold (which we refer to as the horizontal direction),while $\zeta$ is aligned with the radial direction (referred to as the vertical direction) and is the direction along which the HEVI time-integrators are constructed.

 However, the metric terms defined in \eqref{eq:curl-invariant-metrics} and \eqref{eq:jacobian-metrics} are uniquely defined for each element $(e)$ and DOF $(i,j,k)$ which implies that along a vertical column that shares two or more element columns (on the cubed-sphere, this can be 1, 2, 3, or 4 but much higher on an icosahedral grid) the metric terms as defined are $C^0$ continuous. For two of the three methods/storage schemes used in NUMA, this poses no problem. However, for a gridpoint-based continuous Galerkin method (what we call CGc in \cite{abdi:2016,giraldo:2020}) we need to modify the metric terms. This procedure requires the following steps: (i) expand a given metric term $\mathcal{M}$ in tensor-product Lagrange polynomials, (ii) evaluate the weak form within an element, (iii) apply DSS, and (iv) weight the result by the inverse mass matrix.
This procedure ensures that the metric terms are continuous across elements.
 Steps (i) through (iv) outlined above can be defined concisely by the following $L^2$ projection
\be
\label{eq:jacobian-metrics-dss}
\mathcal{M}_I = M^{-1} \bigwedge_{e=1}^{N_e} \inte h_i(\xi) h_j(\eta) h_k (\zeta)  \left( \sum_{l=0}^{N_{\xi}} \sum_{m=0}^{N_{\eta}} \sum_{n=0}^{N_{\zeta}}  h_l(\xi) h_m(\eta) h_n (\zeta) \mathcal{M}^{(e)}_{l,m,n} \right) d\Omega_e
\ee
where $h$ are the 1D basis functions previously defined, $\mathcal{M}_I$ represents the metric terms at the global gridpoint $I$, $\mathcal{M}^{(e)}_{i,j,k}$ represents the metric terms for the element $e$ at the local gridpoint $i,j,k$; e.g., $\mathcal{M}^{(e)}_{i,j,k}$ could represent the components of $\nabla \zeta^{(e)}_{i,j,k}$.
Once we apply this operation to $J$ and $\nabla \zeta$, each vertical column can be integrated independently without the need to DSS along the $(\xi, \eta)$ directions; however, DSS along the $\zeta$ direction is still required but this poses no issue if we partition the grid such that all points of a column are on the same memory space (either an MPI rank or GPU card). 
Note that this process is a one-time cost performed at initialization (or whenever the grid is modified as in dynamically adaptive mesh refinement), so this $L^2$ projection does not negatively impact performance.

\subsection{Testing the Metric Terms}
\begin{figure}
\begin{center}
\begin{subfigure}{0.4\textwidth}
\includegraphics[width=\textwidth]{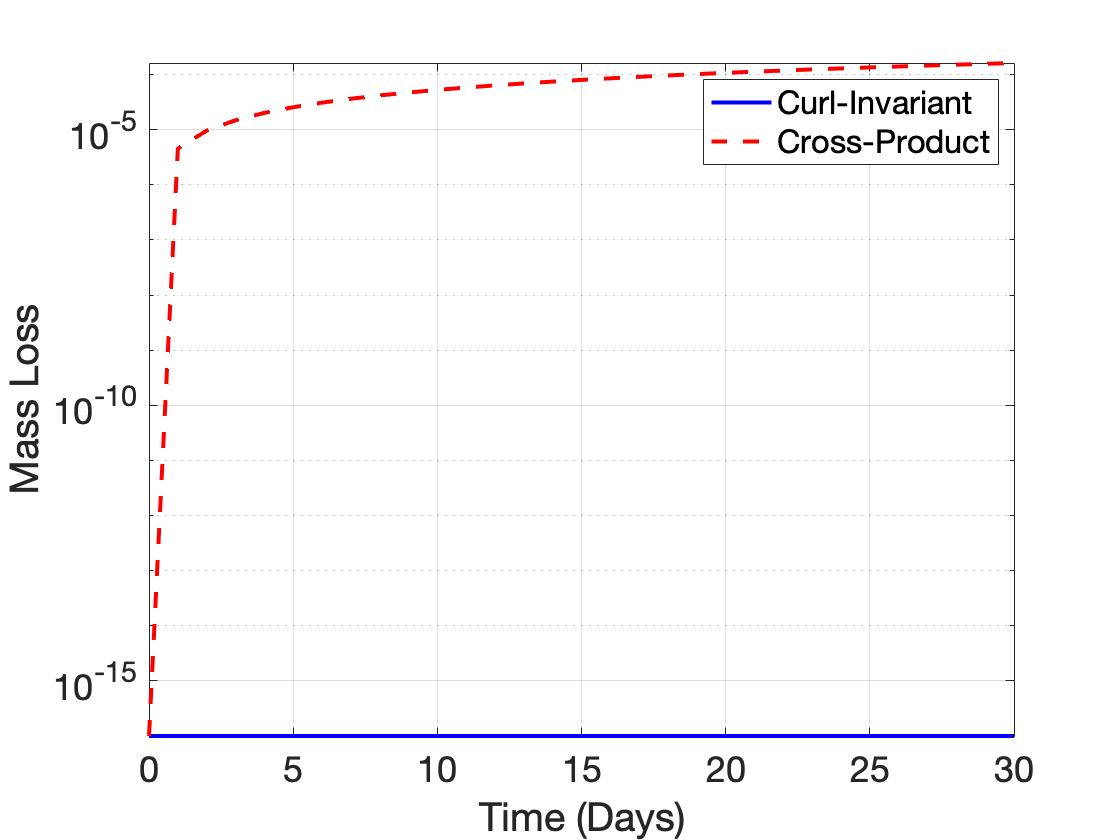}
 \caption{Mountain Case}
 \label{fig:space/mass_conservation/mountain}
\end{subfigure}
\begin{subfigure}{0.4\textwidth}
       \includegraphics[width=\textwidth]{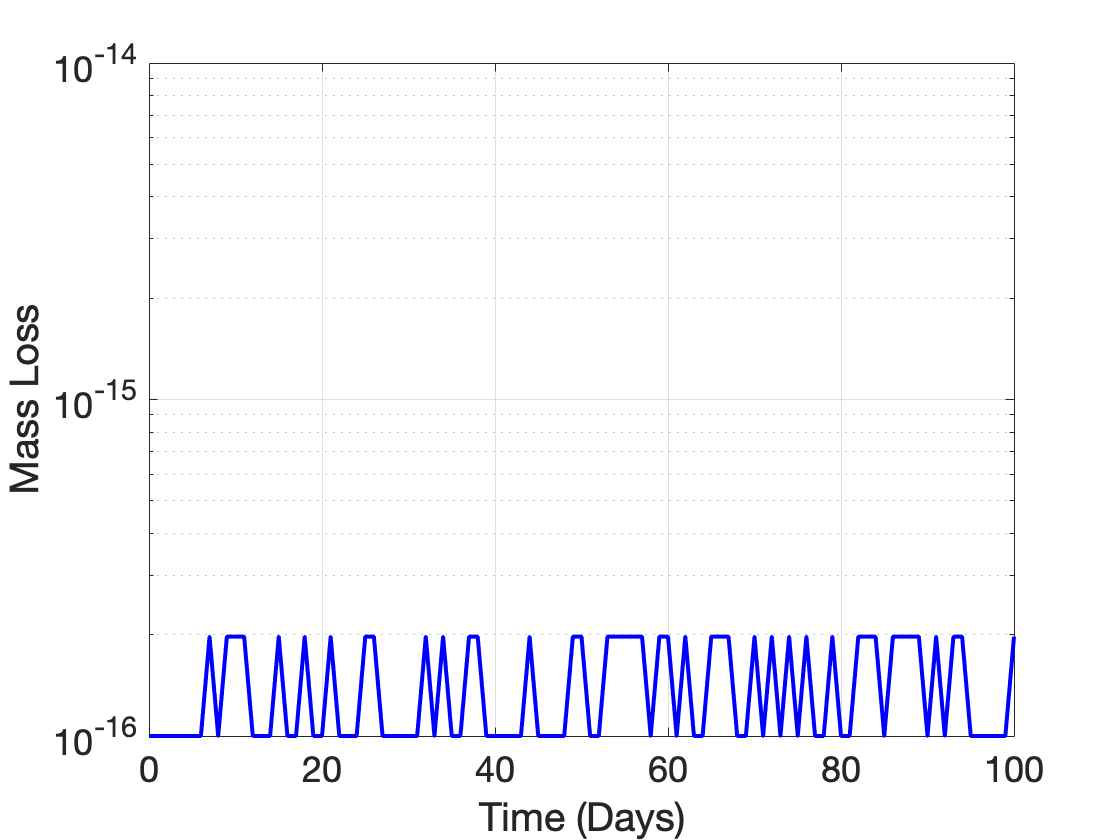}
    \caption{Baroclinic Instability}
    \label{fig:space/mass_conservation/BIS}
\end{subfigure}
\label{fig:space/mass_conservation}
\caption{Mass loss for set 2NC for the (a) static balanced mountain for a 30-day simulation comparing the curl-invariant and cross-product metrics and (b) baroclinic instability for a 100-day simulation using the curl-invariant metrics.}
\end{center}
\end{figure}

Mass conservation is a fundamental property of a dynamical core and important for accurate forecasts, such a surface pressure prediction \cite{thuburn:2008}.  Properly constructing metric terms is a prerequisite for mass conservation. 
To evaluate the curl-invariant metric terms, let us use a test case similar to the mountain case from \cite{tomita:2004}, and compare them with the cross-product metrics. We modify the size of the mountain and apply no initial wind velocity. We define the profile of the axisymmetric mountain as follows
\be
  \mathrm{height} = \frac{h_m}{1+ \left( \frac{r}{a_m} \right)^2}
\ee
 where $r = r_{e} \cos^{-1} \left( \cos \phi \cos \lambda \right)$, $r_e$ is the radius of the earth, with $h_{m}=15$ km and $a_m=5000$ km, with the model top at $z_T = 30$ km. 
 The profile of the mountain has been made artificially large (almost twice the height of Mount Everest) in order to stress the metric terms. 
 For the mountain case, we use $5 \times 5$ elements on each panel of the cubed-sphere grid \cite{ronchi:1996} along with 6 elements in the vertical using polynomial degree $N=4$ in all directions which results in an average grid resolution of 500 km in the horizontal and 1.25 km in the vertical where a terrain-following coordinate is used \cite{galchen:1975}.
 The grid is purposely coarse to test the performance of the metric terms in such resolution regimes.
Figure  \ref{fig:space/mass_conservation/mountain} shows the mass loss for set 2NC for a 30-day simulation with data at a daily cadence; the mass loss is defined as
\be
\label{eq:space/mass_loss}
\mathrm{Mass \, Loss} =\frac{ \abs{M(t) - M(0) }} { M(0) } 
\ee
where
\be
M(t)= \int_{\Omega} \rho(t) d\Omega 
\label{eq:space/mass_integral}
\ee
for any time $t$ integrating along the global domain $\Omega$. The dashed red line in Fig.\ \ref{fig:space/mass_conservation/mountain} representing the cross-product metrics \eqref{eq:cross-product-metrics} hovers near $\order{10^{-5}}$ mass loss. In contrast, the solid blue line for the curl-invariant metrics \eqref{eq:curl-invariant-metrics}, evaluated with \eqref{eq:jacobian-metrics-dss}, conserves mass to machine double precision;  the results for sets 2C and 3C are similar, the only difference is that 3C also conserves total energy. It is for this reason that from here on we use the curl-invariant metrics exclusively. Figure \ref{fig:space/mass_conservation/BIS} shows the mass loss of the curl-invariant metrics for the baroclinic instability (presented in Sec.\ \ref{sec:baroclinic_instability}) for a 100-day simulation using a grid resolution of $104 \times 1.25$ km, where machine double-precision is achieved for the entirety of the 100-day simulation.

\subsection{Tensor-based Hyper-Diffusion}

We introduce numerical  diffusion/hyper-diffusion by adding the following term to the right-hand side of the continuous governing equations
\begin{equation}\label{eq:hyperdiffusion}
\diff{\qvector}{t} = S(\qvector) +  (-1)^{\alpha+1}(\nabla\cdot\vc{\tau}\nabla)^\alpha \qvector,
\end{equation}
where $S(\qvector)$ represents the spatial operators of the compressible Euler equations, $\qvector$ the state variable, $\vc{\tau}$  the viscosity tensor, and $\alpha$ is the order of the hyper-diffusion operator (e.g., $\alpha=1$ is second order standard diffusion and $\alpha=2$ is fourth order hyper-diffusion). In this formulation, we adopt the tensor form of viscosity proposed in \cite{guba:2014b}, but modify the original formulation by multiplying the viscosity tensor $\vc{\tau}$ by each application of Laplacian operator. Applying an EBG discretization to \eqref{eq:hyperdiffusion} yields the semi-discrete form
\begin{equation} 
\label{eq:hyperdiffusion/pde}
  {M} \difft{\qvector}{t} = \mathcal{S}(\qvector) + \mathcal{H}_{\nu}(\qvector),
\end{equation}
with 
\be
\label{eq:hyperdiffusion/operator}
\mathcal{H}_{\nu}(\qvector) = (-1)^{\alpha+1} \qty(\vc{L}_{\nu})^\alpha\qvector,
\ee
where ${M}$ is the global mass matrix, $\mathcal{H}_{\nu}$  is the hyper-diffusion operator, $\vc{L}_{\nu}$ is the Laplacian matrix with viscosity, and $\mathcal{S}$ is the discrete EBG representation of $S$. Repeated application of the Laplacian operator gives rise to higher even-order derivatives \cite{giraldo:1999,kim:2008,giraldo:2020}. The Laplacian matrix $\vc{L}_{\nu}$ is calculated using the contravariant form in terms of the reference coordinates $\vc{\xi}$ on the element $\Omega_e$ as follows
\be
\label{eq:hyperdiffusion/Laplacian}
  \vc{L}_{\nu} = \int_{\Omega_e} \nabla\psi_i \cdot \vc{\tau}\nabla \psi_j \, d\Omega_e
\ee
where $\nabla$ is defined by \eqref{eq:grad_x}. Since $\nabla \psi=\vc{J}^{-\mathcal{T}} \nabla_{\xi} \psi$ with $\nabla_{\xi}$ defined by \eqref{eq:grad_ksi}, $\vc{J}=\pdv{\vc{x}}{\vc{\xi}}$, and $\vc{J}^{-\mathcal{T}} = \left( \vc{J}^{\mathcal{T}} \right)^{-1}$, we rewrite \eqref{eq:hyperdiffusion/Laplacian} as 
\begin{align}
  \vc{L}_{\nu}
&= \int_{\wh{\Omega}} \left( \nabla_{\xi}\psi_i \right) \transpose (\vc{J}^{-1} \vc{\tau}\vc{J}^{-\mathcal{T}}) \nabla_{\xi}\psi_j \, J \, d\wh{\Omega} \notag \\
&= \int_{\wh{\Omega}} \left( \nabla_{\xi}\psi_i \right) \transpose \vc{G}_{\nu} \nabla_{\xi}\psi_j \, J \, d\wh{\Omega},
\end{align}
where $\wh{\Omega}$ denotes the reference element, $J$ the metric Jacobian defined by \eqref{eq:jacobian-metrics}, and the viscous metric tensor $\vc{G}_{\nu}$ is defined as
\begin{equation}\label{eq:viscous_metric}
  \vc{G}_{\nu}=\vc{J}^{-1} \vc{\tau} \vc{J}^{-\mathcal{T}}. 
\end{equation}
Formulating hyper-diffusion using \eqref{eq:hyperdiffusion/pde},  \eqref{eq:hyperdiffusion/operator}, and \eqref{eq:hyperdiffusion/Laplacian} only requires storing a single metric tensor $\vc{G}_\nu$ instead of both a viscosity tensor $\vc{\tau}$ and the viscous metric tensors, thereby reducing memory consumption. 
We employ the following eigendecomposition of the viscosity tensor $\vc{\tau}$ from \cite{guba:2014b}
\begin{equation}\label{eq:tau_tensor}
\vc{\tau}=\vc{JE} \begin{pmatrix}
\nu_1\lambda_1^{-1} & 0 & 0 \\ 0 & \nu_2\lambda_2^{-1} & 0 \\ 0 & 0 & \nu_3\lambda_3^{-1}
\end{pmatrix}(\vc{JE})^{\mathcal{T}},
\end{equation}
where $\left( \nu_1,\nu_2,\nu_3 \right)$ are the viscosity parameters that act in the principal axes of the metric tensor $\vc{G}=\vc{J}^{-1}\vc{J}^{-\mathcal{T}}$, and $\left( \lambda_1^{-1},\lambda_2^{-1},\lambda_3^{-1} \right)$ are the eigenvalues of $\vc{G}$ in ascending order ($\lambda_1>\lambda_2>\lambda_3$). The matrix $\vc{E}$ contains the normalized eigenvectors $\left( \vc{e}_1, \vc{e}_2, \vc{e}_3 \right)$ of $\vc{G}$ in its columns, i.e., $\vc{E}=\left[\vc{e}_1 | \vc{e}_2 | \vc{e}_3\right]$. $\vc{E}$ is orthonormal because $\vc{G}$ is symmetric. The square root $\sqrt{\lambda}$ can be interpreted as the stretching from the reference to the physical element, i.e., $\Delta x=\sqrt{\lambda}\Delta\xi=2\sqrt{\lambda}$ for $\xi\in[-1,1]$, which implies that $\sqrt{\lambda}$ characterizes the length scale of an element.

Substituting \eqref{eq:tau_tensor} into \eqref{eq:viscous_metric} yields the eigendecomposition
\begin{align}
\label{eq:viscous_metric_eigen}
\vc{G}_{\nu} &= \vc{E} \begin{pmatrix}
\nu_1\lambda_1^{-1} & 0 & 0 \\ 0 & \nu_2\lambda_2^{-1} & 0 \\ 0 & 0 & \nu_3\lambda_3^{-1}
\end{pmatrix} \vc{E}^{\mathcal{T}} \notag \\
&= \nu_1\lambda_1^{-1} (\vc{e_1}\otimes\vc{e_1}) + \nu_2\lambda_2^{-1} (\vc{e_2}\otimes\vc{e_2}) + \nu_3\lambda_3^{-1} (\vc{e_3}\otimes\vc{e_3}).
\end{align}
The viscosity parameters $\nu_i$ have physical dimensions $L^2/T^{1/\alpha}$, while kinematic (physical) hyperviscosity $\nu_p$  has dimensions $L^{2\alpha}/T$.
The viscosity parameters can be scaled by the element size $\Delta x$ and time-step $\Delta t$ as follows
\be
\label{eq:viscosity_scaled}
\nu_i \equiv (c_i)^{1/\alpha}\frac{(\Delta x_i)^2}{(\Delta t)^{1/\alpha}} = (c_i)^{1/\alpha}\frac{ 4 \lambda_i}{N^2 (\Delta t)^{1/\alpha}}, \, i=1,\ldots,3
\ee
where $c_i=\left( c_1,c_2,c_3 \right)$ are the dimensionless coefficients that independently control the amount of numerical hyper-diffusion in the principal directions and $\Delta x_i=2\sqrt{\lambda_i}/N$ is the characteristic length for an element of polynomial order $N$ (assumed constant to simplify the discussion); to recover the physical viscosity $\nu_p$ we need to exponentiate \eqref{eq:viscosity_scaled} by $\alpha$.

In many atmospheric applications, a 3D mesh is generated using thin elements. These elements are characterized by two significantly large eigenvalues and a much smaller eigenvalue, $\lambda_1\approx\lambda_2\gg \lambda_3$, due to their high aspect ratio. In this case, we treat viscosity anisotropically by setting $\nu_H=\nu_1=\nu_2$ for the horizontal direction and $\nu_V=\nu_3$ for the vertical direction, assuming that $\lambda_H=\lambda_1=\lambda_2$ and $\lambda_V=\lambda_3$. This setting simplifies the expression for $\vc{G}_\nu$ in \eqref{eq:viscous_metric_eigen} as a combination of horizontal and vertical contributions, as such
\begin{equation}
\vc{G}_{\nu} = \nu_H\lambda_H^{-1} ( \vc{e_1}\otimes\vc{e_1} + \vc{e_2}\otimes\vc{e_2}) + \nu_V\lambda_V^{-1} (\vc{e_3}\otimes\vc{e_3}).
\end{equation}

\section{Horizontally Explicit Vertically Implicit (HEVI)}
\label{sec:hevi}

With the contravariant form of the governing equations and spatial discretization both defined, we are now in a position to describe the HEVI time-integration strategy.
We base the HEVI approach on Additive Runge-Kutta (ARK) methods, although other time-integration approaches could also be used (e.g., multi-step methods). We describe three forms of HEVI methods: two nonlinear (NHEVI-GMRES and NHEVI-LU) methods and one linear (LHEVI) method.  
The difference between linear and nonlinear HEVI is that in nonlinear HEVI (NHEVI), all the vertical terms are solved nonlinearly in an implicit fashion, as opposed to linearly in linear HEVI (LHEVI); this is described in detail below.

\subsection{Implicit-Explicit (IMEX) Runge-Kutta Method}
Since our HEVI method is based on additive Runge-Kutta methods, let us first describe the IMEX Runge-Kutta method.  Consider solving the system of ordinary differential equations (ODEs)
\begin{equation}\label{eq:imex}
  \difft{\vc{q}}{t} = \vc{\Efunction} \left(\vc{q}\right) + \vc{\Ifunction} \left(\vc{q}\right)
\end{equation}
where $\vc{\Efunction}$ denotes the operators that are handled explicitly and $\vc{\Ifunction}$ those that are handled implicitly.

For the state vector $\vc{q}$, the stage values $\vc{Q}_i$ are defined as follows
\begin{equation}
\begin{split}
  \vc{Q}_i = &\vc{q}^n + \Delta t\sum_{j=1}^{i-1} a_{ij} \vc{\Efunction} \left(\vc{Q}_j\right) +  \Delta t\sum_{j=1}^{i} \tilde{a}_{ij} \vc{\Ifunction} \left(\vc{Q}_j\right)  ,\\
\end{split}
\end{equation}
for each stage $i=1,...,s$.
The updated solution is obtained in the following manner
\begin{equation}
  \vc{q}^{n+1} = \vc{q}^n + \Delta t \sum_{i=1}^{s} \left[ b_i \vc{\Efunction} \left(\vc{Q}_i\right) + \tilde{b}_{i} \vc{\Ifunction} \left(\vc{Q}_i\right) \right]
\end{equation}
where \(a_{ij}, b_{i}\) and \(\tilde{a}_{ij}, \tilde{b}_{i}\) are the Butcher tableaux coefficients for the explicit and implicit terms, respectively, and the superscripts $n$ and $n+1$ denote the values at the current and next time-level. In this work, we only consider IMEX methods such that $\mathbf{A}$ is strictly lower triangular and $\tilde{\mathbf{A}}$ is lower triangular, which makes them diagonally-implicit Runge-Kutta (DIRK) methods; to be more precise, we only consider DIRK methods that have the same (singly) diagonal implicit coefficient $ \tilde{a}_{ii}$=constant (SDIRK) and whereby the first stage is explicit (ESDIRK).  The Butcher tableaux are given in Table \ref{table:imex_butcher} for a general Runge-Kutta method with $\mathbf{b}=\tilde{\mathbf{b}}$ which is a necessary condition for conserving linear invariants (e.g., see \cite{giraldo:2013}).
\begin{table}[H]
\caption{Butcher tableaux for the horizontal (left) and vertical (right) terms.}
\centering
\begin{tabular}{c|c}
 \ST $\mathbf{c}$ &  $\mathbf{A}$   \\
 \hline
\ST  & $\mathbf{b}^T$  \\
\end{tabular}
\hspace{0.25in}
\begin{tabular}{c|c}
 \ST $\tilde{\mathbf{c}}$ &  $\tilde{\mathbf{A}}$   \\
 \hline
\ST  & $\tilde{\mathbf{b}}^T$  \\
\end{tabular}
\label{table:imex_butcher}
\end{table}

\begin{figure}
\begin{center}
\begin{subfigure}{0.4\textwidth}
    \includegraphics[width=\textwidth]{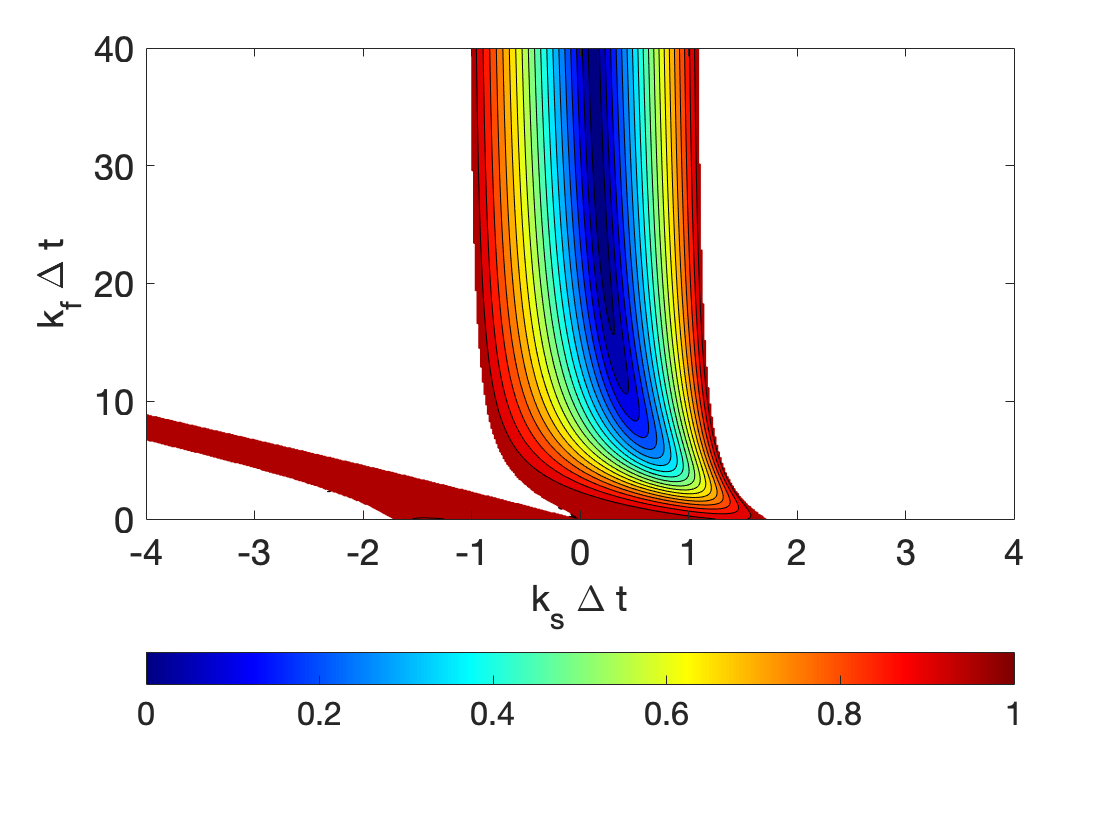}
    \caption{ARK(2,3,2)a}
    \label{fig:time_integration/ark2a_stability}
\end{subfigure}
\begin{subfigure}{0.4\textwidth}
    \includegraphics[width=\textwidth]{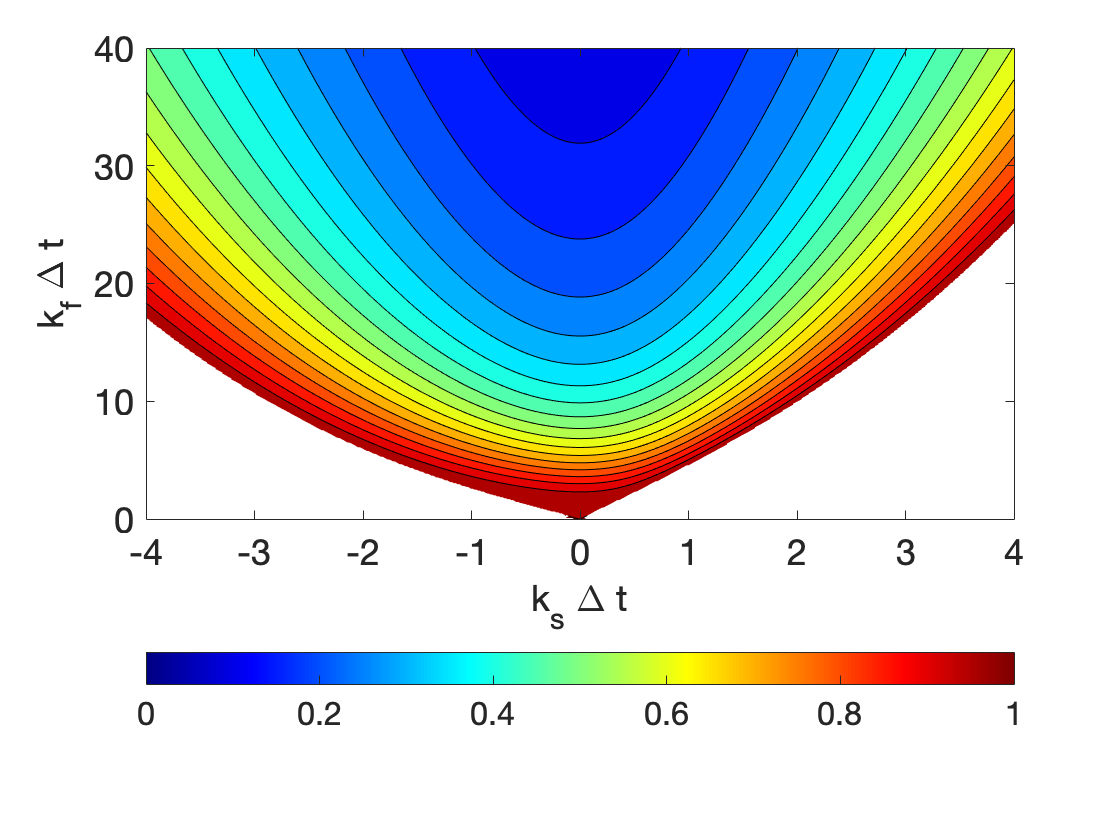}
    \caption{ARK(2,3,2)b}
   \label{fig:time_integration/ark2b_stability}
\end{subfigure}
\caption{IMEX stability region for implicit fast processes ($k_f$) versus explicit slow processes ($k_s$) for the ARK(2,3,2) method. }
\label{fig:time_integration/ark2_stability}
\end{center}
\end{figure}

We explored a range of SDIRK and DIRK methods from the literature including those in \cite{kennedy:2003, giraldo:2013, steyer:2019, vogl:2019, kennedy:2019, guba:2020} and only present those methods that resulted in the most efficient time-to-solution.  To analyze the results of the study, we consider the ordinary differential equation
\be
\difft{q}{t} = \left\{ i k_s q \right \}  + \left[ i k_f q \right]
\ee
where $i=\sqrt{-1}$, $k_s$ and $k_f$ are the frequencies of the slow and fast waves.  If we now analyze the stability region (as we showed in \cite{giraldo:2013}) for implicit fast processes (order $\approx 40$) versus explicit slow processes (order $\approx 4$), for the ARK(2,3,2)\footnote{We use the IMEX designation ARK(i,e,o) from \cite{ascher:1997,pareschi:2005} where i,e,o denote the number of implicit and explicit stages, and the order of the method, respectively.} method presented in \cite{giraldo:2013} we arrive at the stability regions illustrated in Fig.\  \ref{fig:time_integration/ark2_stability}; ARK(2,3,2)a was proposed in \cite{giraldo:2013} ($a_{32}=\frac{1}{6}\left(3 + 2\sqrt{2}\right)$) due to the large explicit stability region along the imaginary axis and for its accuracy (eliminates the second-order error for the explicit components). However, the choice $a_{32}=\frac{1}{2}$ has been the \emph{de facto} method used in NUMA (see \cite{giraldo:2020}, Ch. 20, p. 464) due to its larger allowable time-step. The conclusion from our time-integrator study is that those methods that have the wedge-shape stability region as in Fig.\ \ref{fig:time_integration/ark2b_stability} admit larger time-steps than those methods with the shape shown in Fig.\ \ref{fig:time_integration/ark2a_stability}. It is for this reason that we limit our study to IMEX methods with stability regions similar to those shown in Fig.\ \ref{fig:time_integration/ark2b_stability}.

Below, we describe two types of HEVI methods; in NHEVI 
at each stage value \(\vc{Q}_i\) calculation, a nonlinear implicit system of equations must be solved.  Within NHEVI, we describe two variants: NHEVI-GMRES whereby the linear problem is solved using GMRES and NHEVI-LU where the operator Jacobian (defined in Sec.\ \ref{sec:nhevi-gmres}) is constructed and then solved using direct (LU) solvers.  In contrast, in LHEVI, only a linear system needs to be solved; we can decompose LHEVI into two separate classes: linearization over a time-independent reference state (RS) and linearization over a time-dependent previous solution (PS), although we only focus on LHEVI-PS.  Let us now describe the NHEVI and LHEVI methods.

\subsection{Nonlinear HEVI (NHEVI)}
To motivate the HEVI methods, let us consider a time-integration strategy whereby the stiffness of the PDE is derived from the vertical processes as in the case of global atmospheric simulations where the vertical resolution is much finer than the horizontal. This arises because the spatial scales are so different in these directions where the horizontal dimension is on the order of $2 \pi r_e =  2\pi \cdot 6.371 \times 10^3$ km $\approx 40,000$ km, whereas the vertical domain (in earth weather simulations) is $\approx 40$ km which gives a ratio of 1000:1.  For this reason, the grid resolution in the vertical direction is much finer than in the horizontal. For these types of domains, it behooves us to treat the vertical direction differently and so let us consider treating it implicitly and fully nonlinearly.  
We now describe the details of the NHEVI approach.

Let $\qvector$ be the state vector 
and consider splitting the right-hand side of the PDE into horizontal and vertical terms as follows
\begin{equation}\label{eq:hevi-pde}
  \diff{\vc{q}}{t} = \vc{H} \left(\vc{q}\right) + \vc{V} \left(\vc{q}\right)
\end{equation}
where 
we now seek to solve
\begin{equation}
\begin{split}
  \vc{Q}_i = &\vc{q}^n + \Delta t\sum_{j=1}^{i-1} a_{ij} \vc{H}  \left(\vc{Q}_j\right) +  \Delta t\sum_{j=1}^{i} \tilde{a}_{ij} \vc{V} \left(\vc{Q}_j\right)  ,\\
\end{split}
\end{equation}
where $\vc{H}$ denotes the horizontal terms that are handled explicitly and $\vc{V}$  the vertical terms that are handled implicitly.
Even though we use Cartesian coordinates to represent the velocity field, we are able to decompose the horizontal and vertical directions; we make use of the transformations described in Sec.\ \ref{sec:contravariant_form}.
Before continuing, we need to define the continuous horizontal and vertical operators.

\subsubsection{HEVI Horizontal and Vertical Operators}
In order to define the necessary horizontal and vertical operators for set 2NC, we first expand the divergence operator in the continuity equation \eqref{eq:set2nc_mass} using the contravariant formulation given by \eqref{eq:div_contra} and the gradient operators in \eqref{eq:set2nc_momentum} and \eqref{eq:set2nc_energy} using \eqref{eq:grad}.
The horizontal components given by the $\xi$ and $\eta$ derivatives are collected in 
\begin{equation}
\label{eq:hevi/horizontal_operators}
  \vc{H}(\vc{q}) =  - \begin{pmatrix}
       \frac{1}{J} \nabla_H \cdot \left( J \rho \vc{u}^{H} \right) \\
    \vc{u}^{H} \cdot \nabla_{H} \vc{u} + \frac{1}{\rho} \left(\diff{P}{\xi}\nabla \xi + \diff{P}{\eta}\nabla \eta\right) +  \left(\diff{\Phi}{\xi}\nabla \xi + \diff{\Phi}{\eta}\nabla \eta \right) + 2\omega \wh{\vc{\zeta}} \times \vc{u}  \\
    \vc{u}^{H} \cdot \nabla_{H} \theta
  \end{pmatrix} + \mathcal{H}_{\nu}(\qvector)
\end{equation}
where \(\vc{u}^{H} = \left(u^\xi, u^\eta, 0\right)^\mathcal{T}\), $\nabla_{H} = (\diff{}{\xi}, \diff{}{\eta},0)^\mathcal{T}$, and $\wh{\vc{\zeta}}$ is the unit vector in the ${\vc{\zeta}}$ direction; we include the hyper-diffusion operator in the explicit operator even though it includes both horizontal and vertical components.  The vertical components given by the $\zeta$ derivatives are 
\begin{equation}
\label{eq:hevi/vertical_operators}
  \vc{V}\left(\vc{q}\right) = - \begin{pmatrix}
    \frac{1}{J} \diff{ }{\zeta} \left( J \rho u^\zeta \right) \\
    u^\zeta \diff{\vc{u}}{\zeta} +  \frac{1}{\rho}\diff{P}{\zeta}\nabla\zeta + \diff{\Phi}{\zeta}\nabla\zeta \\
    u^\zeta \diff{\theta}{\zeta}
  \end{pmatrix}.
\end{equation}
Since all the vertical operators, including vertical advection, are included in $\vc{V}$, our scheme is similar to the HEVI-A approach proposed in \cite{gardner:2018}. The horizontal and vertical operators for sets 2C and 3C are defined in \ref{sec:appendix/set2c_set3c/hevi_operators}.

\subsubsection{NHEVI-GMRES}
\label{sec:nhevi-gmres}
In order to compute the stage value \(\vc{Q}_i\), we reformulate the problem similarly to the multidimensional root-finding problem: \( \vc{F} \left(\vc{Q}_i\right) = \vc{0}\), where the nonlinear function \( \vc{F} \) is defined as
\begin{equation}
  \vc{F}\left(\vc{q}\right) = \vc{q} - \Lambda \vc{V} \left(\vc{q}\right) - \vc{R} \ ,
  \label{eq:hevi/nl_function}
\end{equation}
with \(\Lambda =\tilde{a}_{ii}  \Delta t  \) and 
\[
\vc{R} =  \vc{q}^n + \Delta t\sum_{j=1}^{i-1}\left[ a_{ij} \vc{H} \left(\vc{Q}_j\right) + \tilde{a}_{ij} \vc{V} \left(\vc{Q}_j\right)\right].
\]

Given an initial guess \( \vc{q}_0 \), we use the Newton-Raphson method to generate a better approximation to the solution through successive iterations $m$
\begin{equation}
\label{eq:Newton-system-update}
  \vc{q}_{m+1} = \vc{q}_m +\vc{\delta q}_m
\end{equation}
where \(\vc{\delta q}_m\) is computed by solving the linear system
\begin{equation}
\label{eq:Newton-system}
\vc{ \mathbb{J}}( \vc{q}_m) \, \vc{\delta q}_m = - \vc{F} \left(\vc{q}_m\right)
\end{equation}
with
\[
\vc{\mathbb{J}}( \vc{q}_m) = \diff{\vc{F}}{\vc{q}} (\vc{q}_m)
\]
defining the operator Jacobian of the function \( \vc{F} \) evaluated at \(\vc{q}_m\).
In order to avoid building the operator Jacobian matrix (see, e.g., \cite{reisner:2003,knoll:2004}), we use an approximation for the action of the operator Jacobian on the vector \(\vc{\delta q}\) using a Fr\'{e}chet derivative, defined as
\begin{equation}
\label{eq:Frechet-derivative}
 \vc{\mathbb{J}}( \vc{q}_m) \, \vc{\delta q}_m \approx \frac{  \vc{F}\left(\vc{q}_m + \vc{\mathcal{E}} \vc{\delta q}_m \right) - \vc{F} \left(\vc{q}_m \right) }{\vc{\mathcal{E}}}
\end{equation}
where \(\vc{\mathcal{E}}\) is some small perturbation.
\(\vc{\mathcal{E}}\) can be a small scalar value, relative to the norm of \(\vc{q}_m\), but we use different magnitudes for each of the state variables $\qvector$ as described in \cite{vogl:2019}. In this case, \(\vc{\mathcal{E}}\) is defined as a diagonal matrix.

Substituting \eqref{eq:Frechet-derivative} into \eqref{eq:Newton-system} and rearranging yields the final form for the Jacobian-free Newton-Krylov (JFNK) method
\begin{equation}
\label{eq:Newton-system-final}
 \vc{res} \equiv \frac{  \vc{F} \left(\vc{q}_m + \vc{\mathcal{E}} \vc{\delta q}_m \right)  + \left( \vc{\mathcal{E}} - 1 \right) \vc{F} \left(\vc{q}_m\right) }{\vc{\mathcal{E}}} \approx \vc{0}
\end{equation}
that is written in residual form 
that is suitable for the application of the generalized minimal residual (GMRES) \cite{saad:1996} method.  The evaluation of the residual in \eqref{eq:Newton-system-final} only requires one evaluation of $\vc{F}$, which has a similar cost of evaluating the right-hand-side (RHS) in an explicit time-step; we say one and not two because $\vc{F}\left(\vc{q}_m\right)$ is already known. This is a key concept to the power and popularity of JFNK methods.
Finally, since only the vertical terms are solved implicitly, we can solve the system for each column of the mesh independently, which avoids unnecessary off-process data communication and data storage. Because this approach uses a nonlinear HEVI method with GMRES, we refer to it as NHEVI-GMRES.

For set 2NC since the state vector is \( \vc{q} = \left( \rho, \vc{u}^\mathcal{T}, \theta \right)^\mathcal{T} \), then we seek the root of the function
\begin{equation}
\label{eq:hevi5-function}
 \vc{F} \left( \vc{q} \right) = 
 \begin{pmatrix}
  \rho + \Lambda \left(  \frac{1}{J} \diff{\left(J \rho u^\zeta\right)}{\zeta} \right) + R_\rho \\
  \vc{u} + \Lambda \left(u^\zeta\diff{\vc{u}}{\zeta} + \frac{1}{\rho}\diff{P}{\zeta}\nabla\zeta + \diff{\Phi}{\zeta}\nabla\zeta \right) + \vc{R}_u \\
  \theta + \Lambda \left( u^\zeta\diff{\theta}{\zeta} \right) + R_\theta
 \end{pmatrix}
\end{equation}
where \(u^\zeta\) is the third component of the contravariant velocity \(\vc{u}^
\xi = \left(u^\xi, u^\eta, u^\zeta\right)\). Recall that we build the grid such that the \(\zeta\) direction of each element is aligned with the vertical component of the grid, regardless of the geometry (e.g., flow on a sphere or in a box).
For sets 2C and 3C, the function $\vc{F}$ is defined by \eqref{eq:hevi/nl_function} with the corresponding state vector $\qvector$ and right-hand side function $\vc{R}$.

The NHEVI-GMRES method  presents subtle issues including a heavy computational cost.  
To reduce the cost of this algorithm we may consider the following strategies.  The number of Newton iterations can be reduced by increasing the value of the stopping criterion, but this strategy is limited since the solution needs to be accurate. The number of GMRES iterations can be reduced by preconditioning which we did not implement.
The strategies that we did pursue include: (1) replacing the GMRES solver in the solution of the linear problem with LU factorization and (2) removing the Newton solver in the nonlinear loop with frequent linearizations; both of these strategies are described in the sections below where we refer to them as NHEVI-LU in Sec.\ \ref{sec:nhevi-lu} and LHEVI in Sec.\ \ref{sec:lhevi}.

\subsubsection{NHEVI-LU}
\label{sec:nhevi-lu}
An alternative to NHEVI-GMRES requires computing the operator Jacobian analytically. In this approach, we rebuild the Jacobian at each Newton step and then solve the resulting matrix problem using an LU factorization; for this reason we refer to it as NHEVI-LU.  

\paragraph{Analytical Operator Jacobian for Set2NC}
The function \( \vc{F} \left(\vc{q}\right)\) for which we need to find the roots in the NHEVI method is given in \eqref{eq:hevi5-function}.
The analytical operator Jacobian of this function, evaluated at a known state vector \( \vc{q}_0 = \left(\rho_0,\vc{u}_0\transpose,\theta_0\right)\transpose \) and applied to the state vector \( \vc{q} = \left(\rho, \vc{u}\transpose, \theta \right)\transpose \), 
is given by
\be
\label{eq:set2nc-jacobian-jmatrix}
\vc{\mathbb{J}}_{2NC}(\vc{q})=\vc{\mathbb{I}} + \Lambda \mathbb{D}^{\zeta} \vc{\mathbb{K}}_{2NC}(\vc{q})
\ee
where 
\[
\mathbb{D}^{\zeta}_{ij} = \left(  - \frac{(\varpi J)_j}{(\varpi J)_i} D^{\zeta}_{ji}, D^{\zeta}_{ij}, D^{\zeta}_{ij}, D^{\zeta}_{ij}, D^{\zeta}_{ij} \right)
\]
 for the 5 prognostic variables
with $\mathbb{I}$ denoting an identity matrix and
{\footnotesize
\begin{equation}
\label{eq:set2nc-jacobian-kmatrix}
 \vc{\mathbb{K}}_{2NC}(\vc{q})= 
 \begin{pmatrix}
 \left( u^\zeta \right)_j  & 
 \left( \rho \zeta_x \right)_j  & 
\left( \rho \zeta_y \right)_j  & 
\left( \rho \zeta_z \right)_j  & 
 0 \\
 - \left( \frac{\zeta_x}{\rho^2}\right)_i P_j + \left( \frac{\zeta_x}{\rho} \right)_i \diff{P_j}{\rho} & 
\left( \zeta_x \right)_i  u_j + u^{\zeta}_i & 
 \left( \zeta_y \right)_i u_j & 
  \left( \zeta_z \right)_i u_j & 
\left( \frac{\zeta_x}{\rho} \right)_i \diff{P_j}{\theta}  \\
- \left( \frac{\zeta_y}{\rho^2}\right)_i P_j + \left( \frac{\zeta_y}{\rho} \right)_i \diff{P_j}{\rho} & 
\left( \zeta_x \right)_i v_j  & 
 \left( \zeta_y \right)_i v_j + u^{\zeta}_i& 
  \left( \zeta_z \right)_i v_j & 
\left( \frac{\zeta_y}{\rho} \right)_i \diff{P_j}{\theta}  \\
- \left( \frac{\zeta_z}{\rho^2}\right)_i P_j + \left( \frac{\zeta_z}{\rho} \right)_i \diff{P_j}{\rho} & 
\left( \zeta_x \right)_i  w_j  & 
 \left( \zeta_y \right)_i w_j & 
  \left( \zeta_z \right)_i w_j + u^{\zeta}_i  & 
\left( \frac{\zeta_z}{\rho} \right)_i \diff{P_j}{\theta}  \\
0  & 
 \left( \zeta_x \right)_i \theta_j & 
  \left( \zeta_y \right)_i \theta_j & 
   \left( \zeta_z \right)_i \theta_j &  
u^{\zeta}_i
 \end{pmatrix}
\end{equation}
}
where $J$ is the metric Jacobian defined in \eqref{eq:jacobian-metrics}, and $(\zeta_x,\zeta_y,\zeta_z)$ are the metric terms along the vertical direction. 
In \eqref{eq:set2nc-jacobian-jmatrix}, we only show the entries at a given gridpoint $i=0,\ldots,N_{\zeta}$ on one element $\Omega_e$ where the use of $D^{\zeta}_{ij}$ represents a strong form derivative while $-D^{\zeta}_{ji}$ denotes a weak form derivative obtained from integration by parts (e.g., see \cite{giraldo:2020}); our experience shows we need the weak form in the continuity equation to formally conserve mass.
Note that since  $\mathbb{D}^{\zeta}$ is defined for all gridpoints $(i,j)=0,\ldots,N_{\zeta}$ then upon applying the differentiation $\mathbb{D}^{\zeta}$ makes these matrices
$\dim(\vc{\mathbb{J}}_{2NC}, \mathbb{I},\vc{\mathbb{K}}_{2NC})=n_{var}\left( N_{\zeta} + 1\right) \times n_{var}\left( N_{\zeta} + 1\right)$. 

Applying numerical integration (multiplying by $\varpi_i J_i$) and DSS along the $\zeta$ direction along all elements $N^{\zeta}_e$, results in the global representation 
$\mathcal{G}\left( \vc{\mathbb{J}}_{2NC} \right)$ defined in the vector space $\mathbb{R}^{M \times M}$ per vertical column (with $6[N_e^{\xi}N_{\xi}+1][N_e^{\eta}N_{\eta}+1]$ columns on a cubed-sphere grid), where $M=n_{var} n_z$, with $n_z=N^{\zeta}_e \left( N_{\zeta}+1 \right)$ denoting the number of gridpoints along a column, and 
$\mathcal{G}$ represents the operator that takes ${\mathbb{J}}_{2NC}$ and constructs its global representation. We do not present the global matrix because its construction will vary depending on the spatial discretization used, however, the operator at a given gridpoint presented in \eqref{eq:set2nc-jacobian-jmatrix} will not (only the definition of $\mathbb{D}^{\zeta}$).

The cost of building 
$\mathcal{G} \left( \vc{\mathbb{J}} \right)_{2NC}$ 
is $\order{n_{var} N_e^{\zeta} (N_{\zeta}+1)^3}$, with similar costs for constructing the Jacobians for  
sets 2C and 3C; these Jacobians are defined in \ref{sec:appendix/set2c_set3c/operator_jacobian}.

\paragraph{Lapack Banded Matrix Routines}
We build the Jacobian matrix directly in LAPACK's band storage matrix format and use subroutines to perform LU factorization and linear system solves that are specialized to banded matrices \cite{anderson:1999}.
In order to use LAPACK's band storage, the user must specify the number of sub-diagonals \(k_l\) and super-diagonals \(k_u\) of the matrix.
For our methods, the number of sub-diagonals is always equal to the number of super-diagonals and is given by
\begin{equation}
 k_u = k_l = n_{var} (N_{\zeta}+1) - 1 \ ,
\end{equation}
where  \(n_{var}=5\) is the number of state variables.

We already saw that the dimension of the sparse representation of $\mathcal{G} \left( \vc{\mathbb{J}} \right)$  for the NHEVI method is a matrix with size $M \times M$; however, the size of its banded storage form is 
$\left[ 3 n_{var} (N_{\zeta}+1) - 2 \right] \times M$
and its entries are related to that of the original global operator Jacobian matrix $\mathcal{G} \left( \vc{\mathbb{J}} \right)$ 
through
\begin{equation}
 \mathcal{G} \left( \vc{ \mathbb{J}}_B \right)_{2k_u+1+i-j, \ j} = \mathcal{G} \left( \vc{\mathbb{J}} \right)_{i, j}
\end{equation}
for \(\mathrm{max}\left(1, \ j-2k_u\right) \leq i \leq \mathrm{min}\left(n_{var}\times n_z, \  j + k_u\right)\). 

\subsection{Linear HEVI (LHEVI)}
\label{sec:lhevi}
An alternative approach to NHEVI is the linear HEVI (LHEVI) method whereby we solve the linear problem 
\begin{equation}
\label{eq:lhevi-pde}
\diff{\vc{q}}{t} = \left\{ \vc{S} \left(\vc{q}\right) - \vc{\Lfunction} \left(\vc{q}\right) \right\}_{EX} +  \left[ \vc{\Lfunction} \left(\vc{q}\right) \right]_{IM}
\end{equation}
where $\vc{S}$ denotes the original right-hand side operator of the PDE while $\vc{\Lfunction}$ represents a linear approximation of the stiff terms; the curly brackets denote the terms that are solved explicitly (EX) whereas the square brackets are those that are solved implicitly (IM). We must consider how to create the linear operator $\vc{\Lfunction}$; this approach is akin to a Rosenbrock \cite{rosenbrock:1963} method since $\vc{\Lfunction}$ is an approximation to the fully implicit nonlinear operator Jacobian $\vc{\mathbb{J}}$.  The simplicity and efficiency of LHEVI can be illustrated by discretizing \eqref{eq:lhevi-pde} in time using a forward/backward Euler method
\begin{align}
\label{eq:lhevi-pde-ark1}
\qvector^{n+1} &= \qvector^n + \Delta t \vc{S} \left(\qvector^n \right) + \Delta t \vc{\Lfunction} \left( \qvector^{n+1} - \qvector^n \right) \\
&=\vc{q}^E + \Delta t \vc{\Lfunction} \left( \qvector^{n+1} - \qvector^n \right), 
\end{align}
subtracting both sides by $\qvector^{n}$
\be
\label{eq:lhevi-pde-ark1-2}
\qvector^{n+1} - \qvector^{n} = \vc{q}^E - \qvector^{n} + \Delta t \vc{\Lfunction} \left( \qvector^{n+1} - \qvector^n \right) ,
\ee
and letting $\qvector_{tt}=\qvector^{n+1} - \qvector^{n}$, $\wh{\qvector}=\vc{q}^E - \qvector^{n}$, and $\Lambda=\Delta t$ yields the final LHEVI form
\be
\label{eq:lhevi-ark-final}
\qvector_{tt} = \wh{\qvector} + \Lambda \vc{\Lfunction} \qvector_{tt}
\ee
which can be rearranged to more easily illustrate the implicit operator as follows
\be
\label{eq:lhevi-ark-final-2}
\left( \vc{\mathbb{I}} - \Lambda \vc{\Lfunction} \right) \qvector_{tt} = \wh{\qvector}.
\ee

The extension to an arbitrary order IMEX additive Runge-Kutta (ARK) method is described in \cite{giraldo:2013} which we summarize now for completeness. We employ \eqref{eq:lhevi-ark-final} with the following definitions 
\begin{align}
\label{eq:lhevi-ark-final-3}
\qvector_{tt} &= \vc{Q}^{(i)} + \sum_{j=1}^{i-1} \frac{\tilde{a}_{ij}-a_{ij}}{\tilde{a}_{ii}} \vc{Q}^{(j)} \notag \\
\wh{\qvector} &= \qvector^E + \sum_{j=1}^{i-1} \frac{\tilde{a}_{ij}-a_{ij}}{\tilde{a}_{ii}} \vc{Q}^{(j)} \\
\qvector^E &= \qvector^n + \Delta t \sum_{j=1}^{i-1} a_{ij} \vc{S} \left( \vc{Q}^{(j)} \right) \notag .
\end{align}

\subsubsection{LHEVI Linear Operator}
To construct a  semi-implicit IMEX method for the nonhydrostatic equations we construct an approximation $\vc{\Lfunction}$ by applying the expansion $\qvector(\vc{x},t)=\qvector_0(\vc{x},t) + \qvector'(\vc{x},t)$ in $\vc{S}$ \cite{baldauf:2021,lee:2021,sridhar:2022,waruszewski:2022,souza:2023}, where $\qvector_0(\vc{x},t)$ represents the solution at a previous time-level $t^n$, and only retaining terms linear in
$\qvector'(\vc{x},t)$.  Doing so, results in the linear operator 
\begin{equation}
    \label{eq:set4nc_1d_imex_PS}
    \vc{\Lfunction}\left(\vc{q}\right) = -\begin{pmatrix}
        \frac{1}{J} \diff{ }{\zeta} \left[ J \left( \rho_0 u^\zeta + \rho u_0^\zeta \right) \right]\\
       u^\zeta_0 \diff{\vc{u}}{\zeta} + u^\zeta \diff{\vc{u}_0}{\zeta} +   \frac{\nabla{\zeta}}{\rho_0} \left[ \diff{ }{\zeta} \left( F_0 \rho + G_0 \theta \right) - \frac{\rho}{\rho_0}\diff{P_0}{\zeta}\right] \\ 
       u^\zeta \diff{ \theta_0}{\zeta} + u_0^\zeta \diff{ \theta}{\zeta} 
    \end{pmatrix}
\end{equation}
which we call LHEVI-PS (previous solution).
A special case of LHEVI-PS arises when we apply a  time-independent linearization defined as $\qvector(\vc{x},t)=\qvector_0(\vc{x}) + \qvector'(\vc{x},t)$; linearizing only about the scalar variables results in the semi-implicit IMEX methods described in \cite{giraldo:2010,giraldo:2013,weller:2013,lock:2014,abdi:2017b,mueller:2018} which we call LHEVI-RS (reference state). The reference state $\qvector_0(\vc{x}) $ is selected to satisfy either a hydrostatic balance and/or geostrophic balance.

In what follows, we build on the NHEVI-LU method described in Sec.\ \ref{sec:nhevi-lu} and use $\vc{\Lfunction}=\mathbb{I} - \vc{\mathbb{J}}$, where  
$\vc{\mathbb{J}}$ is the operator Jacobian of $\vc{F}$ for each equation set; we can in fact see this relationship between $\vc{\Lfunction}$ and $\vc{\mathbb{J}}$ in \eqref{eq:lhevi-ark-final-2}.  We then use LAPACK's banded routines to perform the LU factorization and back substitution as described in Sec.\ \ref{sec:nhevi-lu}. Therefore, in the LHEVI-PS method, we essentially apply the NHEVI-LU method but for only one Newton iteration.
In order to avoid the creation, factorization, and solution of the linear system at every time step, we explore an option to set how often the reference state is modified; this is the reason why we refer to this state as a \emph{previous solution}. This reference state can be any previous state including the current one as in NHEVI.
Unless otherwise specified, we recompute the operator Jacobian every 5 time-steps and show in Sec.\ \ref{sec:results} that some time-integration methods allow for much larger time-steps between operator Jacobian updates (referred to as the update time-step).

\section{Complexity Analysis}
\label{sec:complexity}
The number of vertical levels $n_z$ may need to be increased in an atmospheric model for several applications, such as studying the planetary boundary layer (PBL) or extending a model into the upper atmosphere (e.g., whole-atmosphere modeling) \cite{kelly:2023}.  Hence, knowledge of the computational complexity of the time-integration scheme with respect to $n_z$ is an important consideration.  Let us now describe the complexity of each of the three methods: NHEVI-GMRES, NHEVI-LU, and LHEVI (PS).

\subsection{NHEVI-GMRES}
In NHEVI-GMRES, we must solve the JFNK problem for all $n_{var}=5$ prognostic variables (density, three velocities, and thermodynamic variable) so that $n_z$ gridpoints along a vertical column require the solution of a $5n_z \times 5n_z$ system of nonlinear equations to be solved with complexity
\be
25 n_z^2 N_{Newton}(col) N_{GMRES}(col)^2
\label{eq:nhevi-gmres-complexity}
\ee
per column $col$ per RK stage, where $N_{Newton}$ denotes the number of Newton iterations, and $N_{GMRES}$ the number of GMRES iterations. Because the number of Newton and GMRES iterations vary for each column, the cost per iteration will be that of the most expensive column.  Recall that the number of points along a vertical column is defined as $n_z=N_e^{\zeta} N_{\zeta} + 1$, where $N_{\zeta}$ and $N_e^{\zeta}$ represent the polynomial order and number of elements along the $\zeta$ direction, which is aligned with the radial component of the sphere (recall Fig.\ \ref{fig:space/reference_grid}). 

\subsection{NHEVI-LU}
Next, let us consider reducing the complexity (and computational cost) of NHEVI-GMRES by obviating the need for the GMRES iterations. We achieve this by replacing the iterative solution with a direct one such as LU, although other direct solvers can be substituted.
A standard LU factorization requires $\mathcal{O}\left( 50 n_z^3 \right)$ for the forward reduction and 
$\mathcal{O}\left(25 n_z^2 \right)$  for the back substitution. Hence the complexity of NHEVI-LU is
\be
\left( 50 n_z^3 + 25 n_z^2 \right) N_{Newton}(col)
\label{eq:nhevi-lu-complexity}
\ee
per column per RK stage.
The complexity of this approach can be further reduced by storing the matrix 
in banded form. For a matrix with bandwidth $B$ the per column complexity of the forward-reduction and back-substitution is 
\be
10 n_z B( B + 1)
\label{eq:nhevi-lu-banded-complexity}
\ee
where $B=2 n_{var} (N_{\zeta}  + 1) \approx 10 N_{\zeta} $. 
Therefore, the total per column per RK stage complexity of NHEVI-LU is
\be
100 n_z N_{\zeta} \left( 10 N_{\zeta} + 1 \right) N_{Newton}(col).
\label{eq:nhevi-lu-complexity}
\ee
Comparing the complexity of NHEVI-LU \eqref{eq:nhevi-lu-complexity} with NHEVI-GMRES \eqref{eq:nhevi-gmres-complexity}, we see that the cost of NHEVI-GMRES scales quadradically with $n_z$ and $N_{GMRES}$ whereas NHEVI-LU scales linearly with $n_z$ and does not require GMRES. So, as the number of vertical levels $n_z$ increases, NHEVI-LU requires fewer FLOPS than NHEVI-GMRES, making NHEVI-LU more efficient.

\subsection{LHEVI}
The advantage of LHEVI can be described by reconsidering the total per column complexity of NHEVI-LU, which we modify as follows for LHEVI
\begin{subequations}
\be
1000 n_z  N_{\zeta}^2 
\label{eq:lhevi-complexity-forward}
\ee
\be
100 n_z N_{\zeta}
\label{eq:lhevi-complexity-backsolve}
\ee
\label{eq:lhevi-complexity}
\end{subequations}

\noindent
where \eqref{eq:lhevi-complexity-forward} is the cost of the forward-reduction and is only incurred whenever we update the operator Jacobian, while \eqref{eq:lhevi-complexity-backsolve} is the cost of the back-substitution that is incurred per RK stage. Comparing the NHEVI-LU and LHEVI complexity, we see that LHEVI amortizes the cost of constructing the forward-reduction across all stages and numerous time-steps so that we mostly pay the price of the less computationally expensive back-substitution.

\subsection{Summary of Complexity} 
\begin{figure}
\begin{center}
 \includegraphics[width=0.4\textwidth]{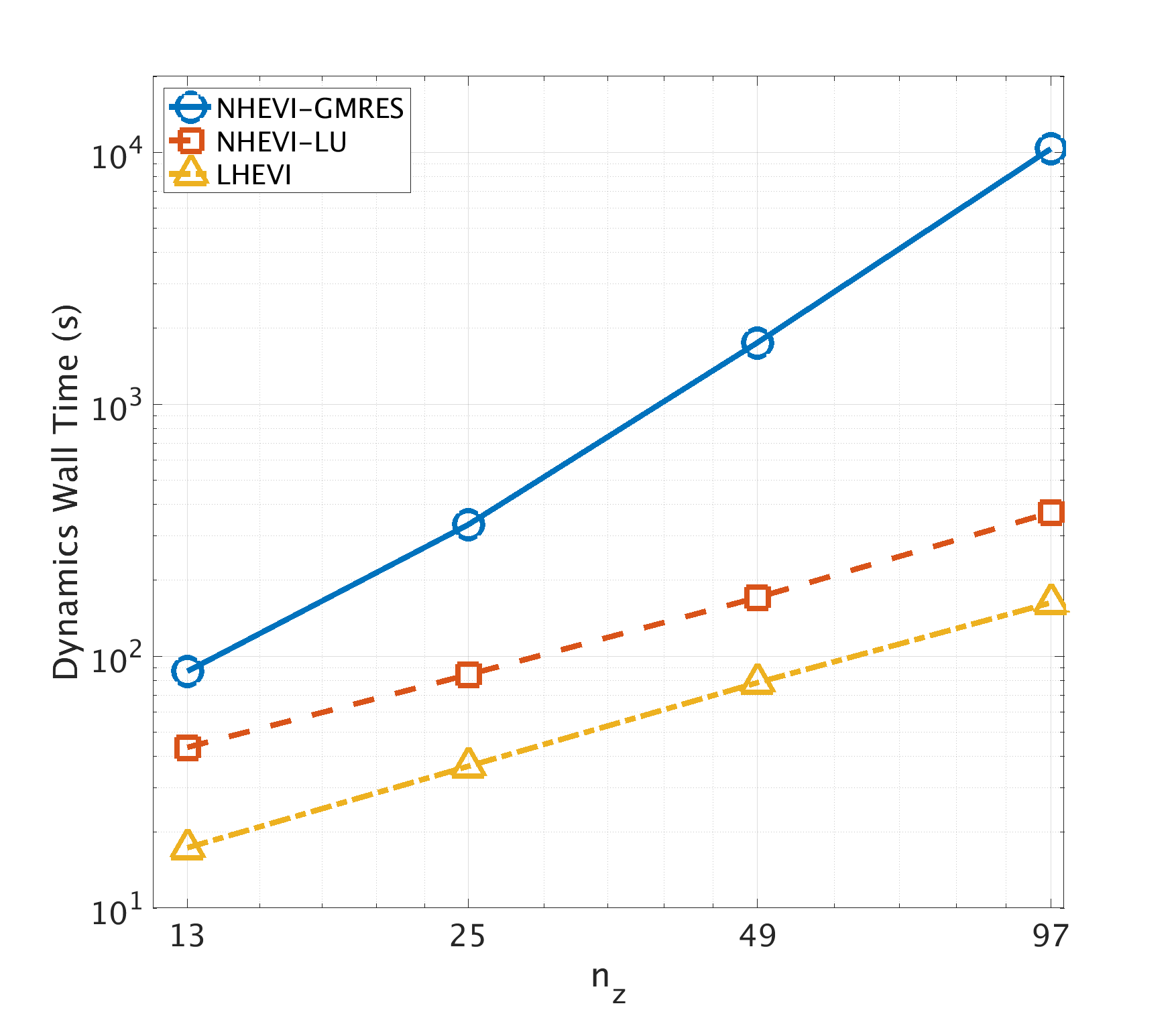}
\caption{Time-to-solution for 0.1 simulation days as a function of vertical gridpoints $n_z$.}
\label{fig:hevi/scaling_nz}
\end{center}
\end{figure}

To summarize, the complexity of both NHEVI-LU and LHEVI scale linearly with $n_z$ while NHEVI-GMRES scales quadratically with $n_z$. To illustrate this behavior, we ran four simulations with $n_z$ ranging from 13 to 97 with each of the three HEVI variants for a short period of time (0.1 days at 104 km horizontal resolution for the baroclinic instability) and only report the simulation time associated with dynamics (and discount initialization time). The results of this study shown in Fig.\ \ref{fig:hevi/scaling_nz}
 clearly show the super-quadratic dependence of $n_z$ for NHEVI-GMRES and the approximately linear dependence for NHEVI-LU and LHEVI. In Sec.\ \ref{sec:results} we discuss a more detailed cost analysis for  the three methods. It should be added that increasing $n_z$ has an adverse effect on the time-to-solution for NHEVI-GMRES because satisfying the stopping criterion on the residual $\vc{res}$ is dependent on the condition number \cite{trefethen:1997}. This increases the computation cost and also the memory consumption, since more Krylov vectors need to be stored.

\section{Time-Integrators}
\label{sec:time_integrators}

In this section we explore five additive Runge-Kutta (ARK) time-integrators taken from the literature covering a range of orders from second through fifth. Our intent here is not to perform an exhaustive study of IMEX methods but rather to 
evaluate a representative sample of time-integrators that gave us the fastest time-to-solution.
The methods that we showcase below are: ARK2 which is the ARK(2,3,2)b method defined in \cite{giraldo:2013,giraldo:2020}, 
ARK3 which is the ARK(3,4,3) method from \cite{kennedy:2003}, 
ARS3 which is the ARK(3,4,3) method from \cite{ascher:1997}, ARK4 which is the ARK(5,6,4) from \cite{kennedy:2003}, and ARK5 which is the ARK(7,8,5) from \cite{kennedy:2019}. Because ARS3 and ARK5 do not possess a dense (or continuous) output formula (the rest of the methods do), we did not use it for constructing an initial guess for the nonlinear solver; this way we perform a fair comparison across all the methods. Dense output only impacts the performance of the NHEVI methods.

\begin{figure}
\begin{center}
\begin{subfigure}{0.19\textwidth}
    \includegraphics[width=\textwidth]{figures/Time_Integrators/ARK2B.png}
    \caption{ARK2}
    \label{fig:time_integration/ark2_stability}
\end{subfigure}
\begin{subfigure}{0.19\textwidth}
    \includegraphics[width=\textwidth]{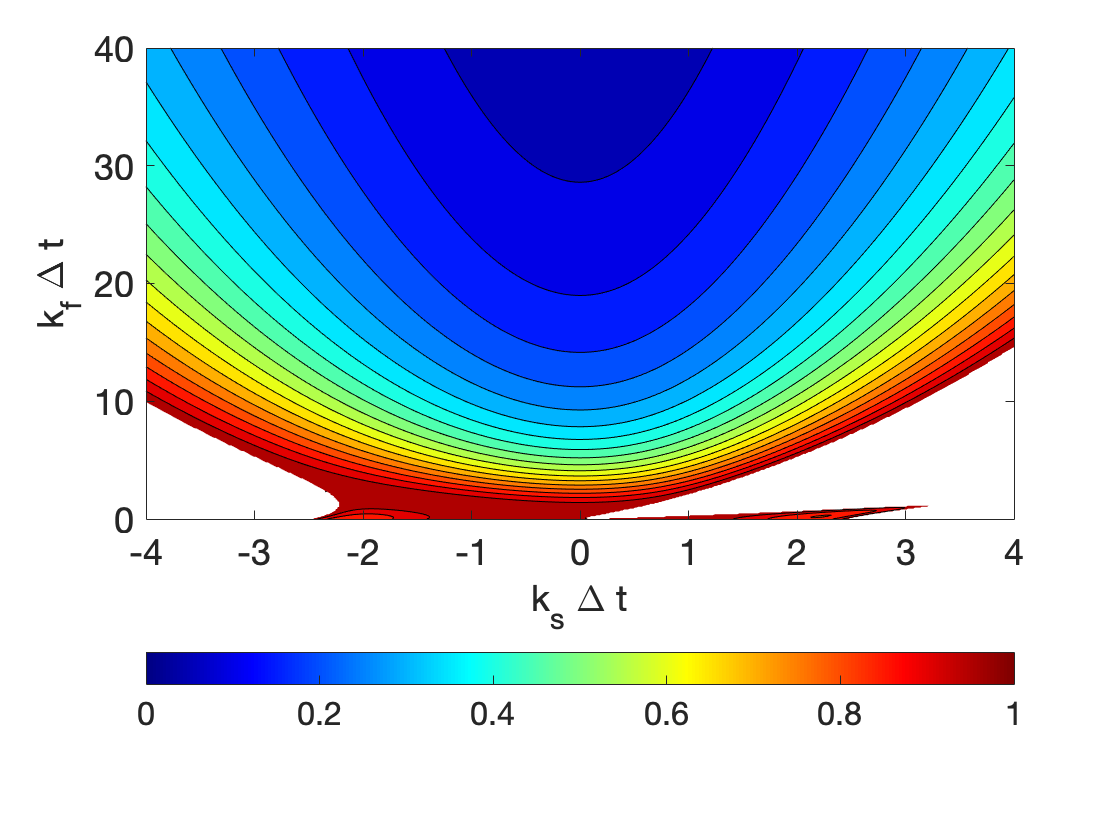}
    \caption{ARK3}
    \label{fig:time_integration/ark3_stability}
\end{subfigure}
\begin{subfigure}{0.19\textwidth}
    \includegraphics[width=\textwidth]{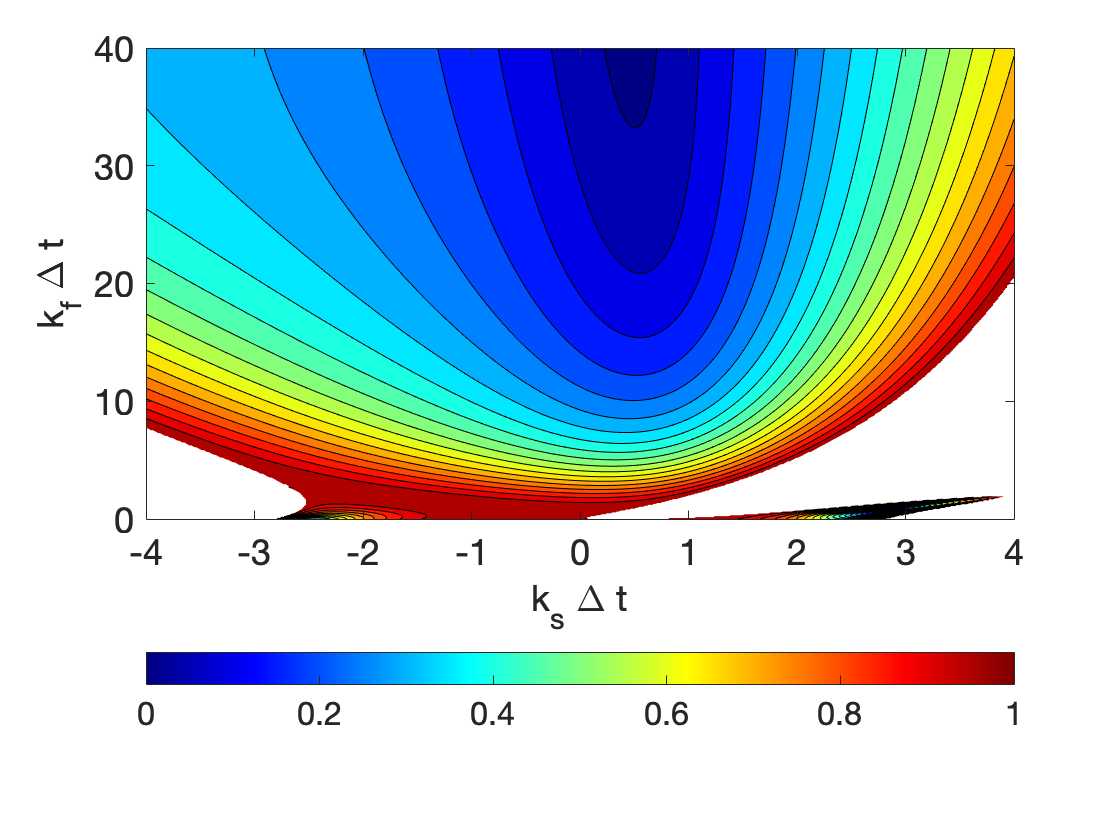}
    \caption{ARS3}
   \label{fig:time_integration/ars343_stability}
\end{subfigure}
\begin{subfigure}{0.19\textwidth}
    \includegraphics[width=\textwidth]{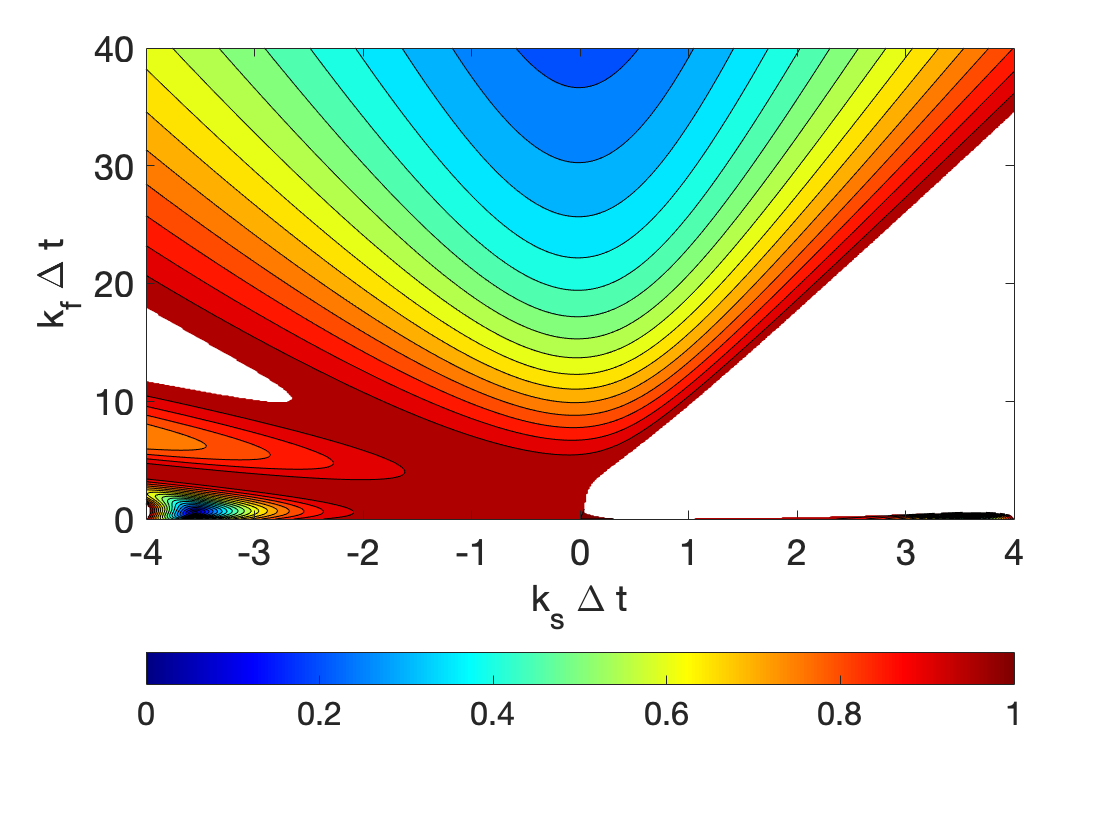}
    \caption{ARK4}
   \label{fig:time_integration/ark4_stability}
\end{subfigure}
\begin{subfigure}{0.19\textwidth}
    \includegraphics[width=\textwidth]{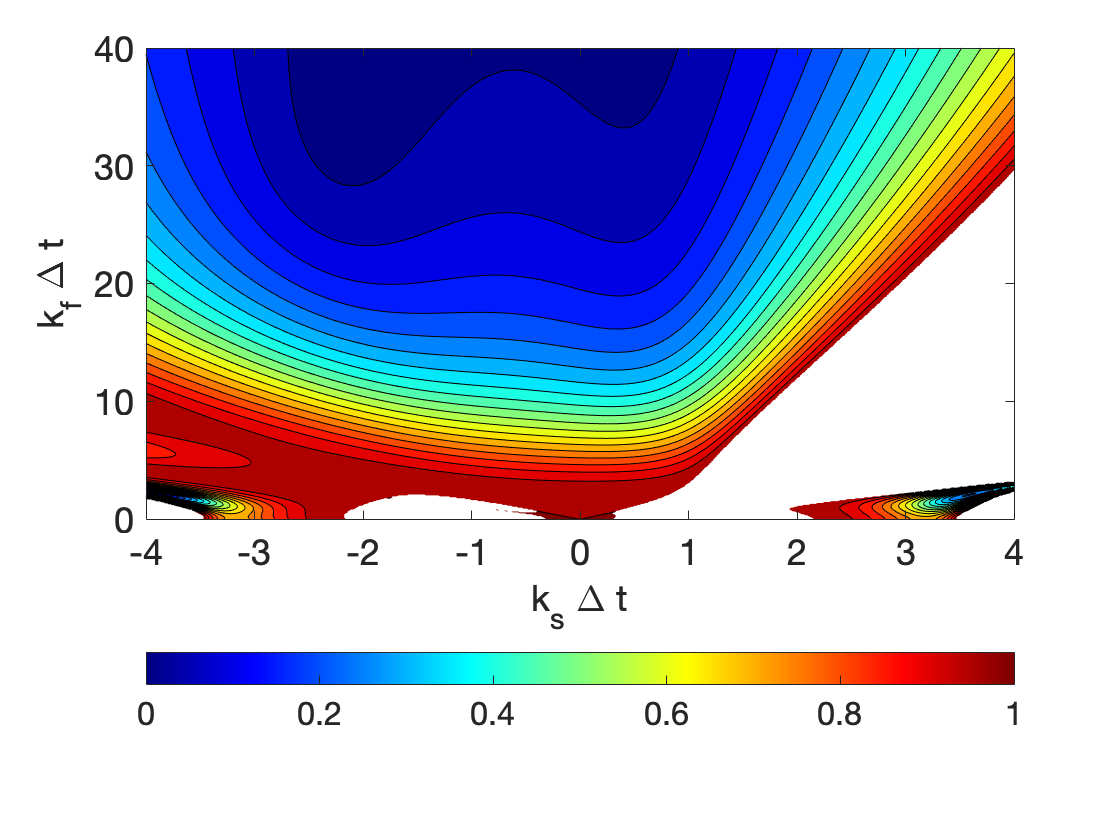}
    \caption{ARK5}
   \label{fig:time_integration/ark5_stability}
\end{subfigure}
\caption{IMEX stability region for implicit fast processes ($k_f$) versus explicit slow processes ($k_s$) for a variety of ESDIRK methods. }
\label{fig:time_integration/esdirk_stability}
\end{center}
\end{figure}

Figure \ref{fig:time_integration/esdirk_stability} shows the stability region of the five methods we consider. Note that they all have the wedge-shape stability region which confers a larger explicit stability region with increasing implicit time-step size; this is the advantage of these types of time-integrators. The disadvantage is that it is difficult to know \emph{a priori} how diffusive the method is with increasing implicit time-step. Comparing Figs.\ \ref{fig:time_integration/ark2a_stability} and \ref{fig:time_integration/ark2b_stability}, we see that for those methods where the implicit time-step does not alter the stability region of the explicit part, we can always be sure to avoid the diffusive part of the region by pushing the time-step to its maximum explicit component (red region in Fig.\ \ref{fig:time_integration/ark2a_stability}). However, this means that these methods will always require a smaller time-step; these are the sorts of trade-offs that need to be considered when choosing a method. In what follows, we aim to use methods that allow for a faster time-to-solution regardless of the dissipation introduced.

Comparing the performance of the various methods using complexity is possible except for the NHEVI-GMRES method. The reason for this is because NHEVI-GMRES requires outer Newton iterations in addition to inner GMRES iterations; the main complication is in the inner GMRES iterations. In the baroclinic instability test below, we require on average 2-7 Newton iterations per time-step per column and $\approx 45$ GMRES iterations per time-step per column (this number includes the iterations for all implicit stages and held constant across all time-integrators).  For NHEVI-LU, computing the cost is more straightforward since the cost of the LU decomposition and backsolve is well-defined; this is also true for LHEVI since there are no Newton iterations.  Therefore, in LHEVI, the cost of the method is directly proportional to the effective time-step of the method which we define as
\[
\Delta t_{eff}=\frac{\Delta t_{max}}{s_{im}}
\]
where $\Delta t_{max}$ is the maximum time-step that the method can use and $s_{im}$ denotes the number of implicit stages.
This means that the number of evaluations of the left- and right-hand side functions are on the order
\[
\mbox{Cost}_{LHEVI}= \frac{ T_{final} }{\Delta t_{eff} }.
\]
To compare whether method 1 is faster than method 2, we compute the ratio
\be
\label{eq:time_integration/R_ratio}
R_{m1\rightarrow m2}= \frac{ \Delta t_{eff, m1} } {\Delta t_{eff, m2} }
\begin{array} {c c}
> 1 & faster \\
< 1 & slower
\end{array}.
\ee
We use these ratios in Sec.\ \ref{sec:results} to compare the cost of the time-integrators.

\section{Numerical Results}
\label{sec:results}

In this section, we present the numerical results for two dry test cases on the sphere.  We begin with the nonhydrostatic baroclinic instability \cite{ullrich:2014} to discern the performance of all the methods and to gauge the maximum usable time-steps for each time-integrator. Finally, we use the inertia-gravity wave problem \cite{tomita:2004}  to compute convergence rates for the time-integrators. 

\subsection{Baroclinic Instability}
\label{sec:baroclinic_instability}
\begin{figure}
\begin{center}
\begin{subfigure}{0.4\textwidth}
     \includegraphics[width=\textwidth]{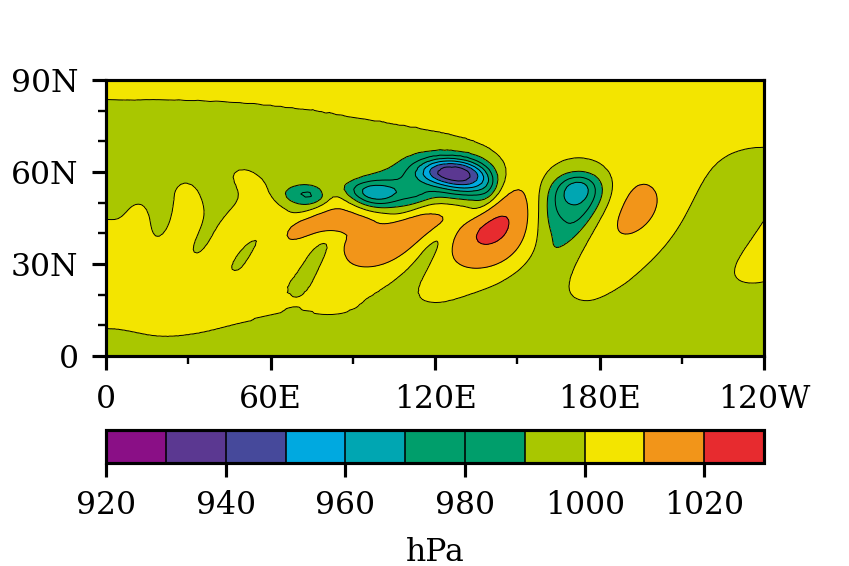}
    \caption{Surface Pressure}
    \label{fig:results/BIS/surface_pressure}
\end{subfigure}
\begin{subfigure}{0.4\textwidth}
       \includegraphics[width=\textwidth]{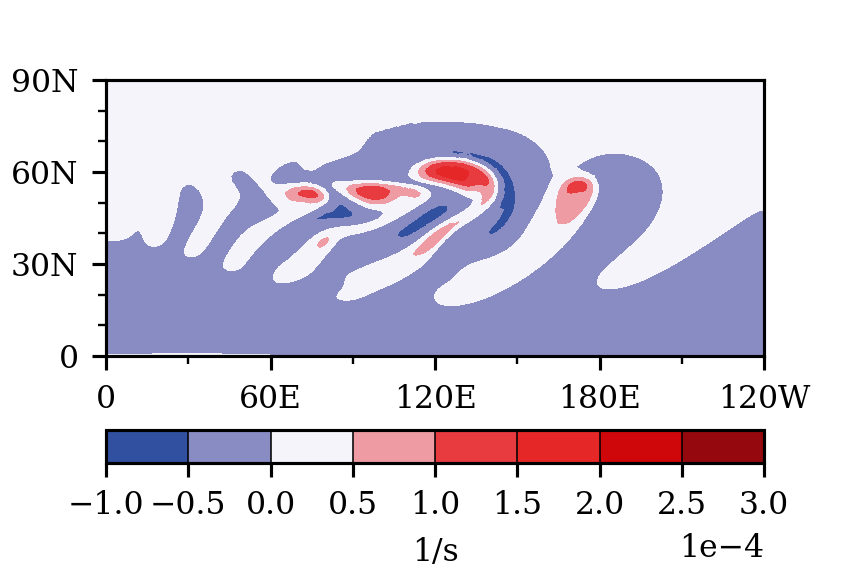}
    \caption{Vertical Vorticity}
    \label{fig:results/BIS/vorticity}
\end{subfigure}
\caption{Baroclinic Instability. Set 2NC with ARK2-LHEVI at grid resolution $104 \times 1.25$ km resolution at 10 days showing: (a) Surface pressure and (b) vertical vorticity at 850 hPa.}
\label{fig:results/BIS}
\end{center}
\end{figure}

The nonhydrostatic baroclinic instability proposed by Ullrich et al.\ \cite{ullrich:2014} is a useful idealized test case for evaluating global, deep-atmosphere, non-hydrostatic dynamical cores.  In particular, this case has been used for previous HEVI studies \cite{gardner:2018,vogl:2019}.  The background state for this case is an exact solution to the deep-atmosphere Euler equations that is in both hydrostatic and geostrophic balance \cite{staniforth:2011}.  A stream-function perturbation is applied, triggers a baroclinic instability, and wave breaking occurs after day 7.  

\begin{figure}
\begin{center}
\begin{subfigure}{0.4\textwidth}
     \includegraphics[width=\textwidth]{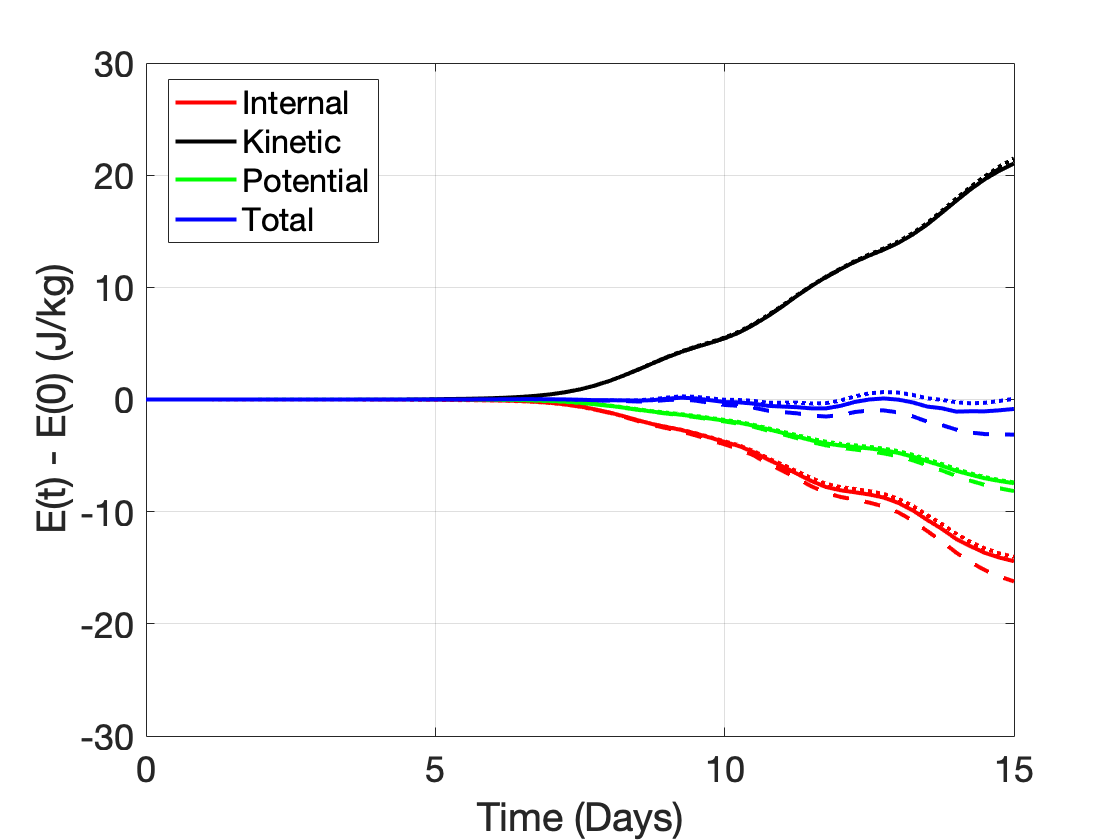}
    \caption{Energy Budget}
    \label{fig:results/BIS/energy_budget}
\end{subfigure}
\begin{subfigure}{0.4\textwidth}
     \includegraphics[width=\textwidth]{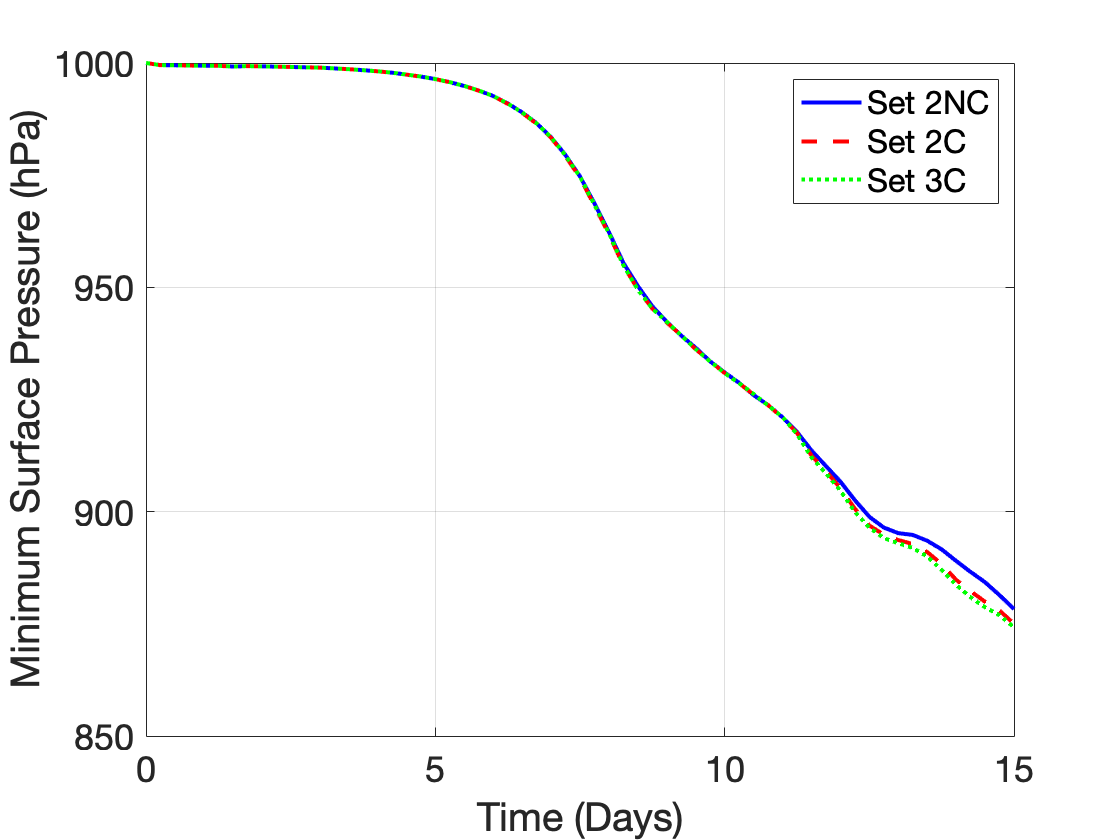}
    \caption{Surface Pressure}
    \label{fig:results/BIS/minimum_pressure}
\end{subfigure}
\caption{Baroclinic Instability. (a) Energy budget (solid line is 2NC, dashed is 2C, and dotted is 3C) and (b) minimum surface pressure for three equation sets. All simulations use ARK2-LHEVI at grid resolution $104 \times 1.25$ km resolution through 15 days.}
\label{fig:results/BIS/energy_pressure}
\end{center}
\end{figure}

Figure \ref{fig:results/BIS} shows contours for the surface pressure and vertical vorticity at 850 hPa for set 2NC using a grid resolution of $104 \times 1.25$ km; the results for the other two equation sets are similar (not shown). The grid resolution corresponds to using polynomial order $N_{\xi}=N_{\eta}=N_{\zeta}=4$ with number of elements $N^{\xi}_e=N^{\eta}_e=24$ and $N^{\zeta}_e=6$.
The surface pressure contours shown in Fig.\ \ref{fig:results/BIS/surface_pressure} and vertical vorticity shown in Fig.\ \ref{fig:results/BIS/vorticity} compare favorably with results from the literature \cite{ullrich:2014,waruszewski:2022,skamarock:2021}. 
 We also show the energy budget and minimum surface pressure for the three equation sets with data at a 6-hour cadence in Fig.\ \ref{fig:results/BIS/energy_pressure}. Figure \ref{fig:results/BIS/energy_budget}
 shows the deviations in internal, kinetic, potential, and total energy from their initial values and are computed in a similar fashion to the mass loss in \eqref{eq:space/mass_loss} and \eqref{eq:space/mass_integral} (without normalizing by the initial values). For sets 2NC and 2C, total energy is constructed from the other three components whereas for set 3C, internal energy is diagnosed (since total energy is a prognostic variable). 
Because the hyperdiffusion applied to the momentum and thermodynamic equations dissipates energy, we do not expect exact conservation for any of the three equation sets; however sets 2NC and, especially, 3C do very well and exhibit low total energy loss. Figure \ref{fig:results/BIS/minimum_pressure} shows that all equation sets yield similar mimimum surface pressure for most of the simulation. The results for both the energy budget and minimum surface pressure compare well with results from the literature (see, e.g., \cite{skamarock:2021}).
With the verification confirmed, let us now turn to the discussion of the performance of the various HEVI methods. 

\begin{figure}
\begin{center}
 \includegraphics[width=0.5\textwidth]{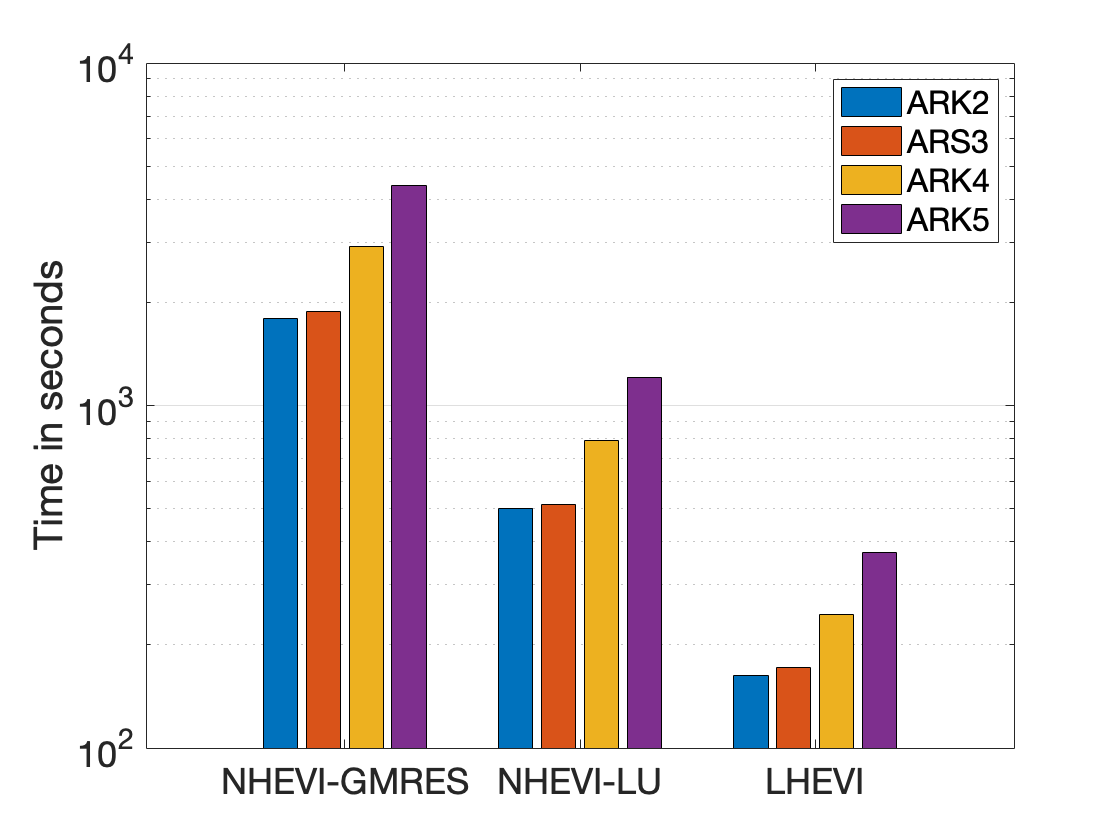}
 \caption{Baroclinic Instability. Time-to-solution (in seconds) for set 2NC for a 10-day simulation at $104 \times 1.25$ km resolution. Each IMEX method uses its maximum usable time-step configuration presented in Table  \ref{table:results/BIS/time-to-solution}.}
 \label{fig:results/BIS/time-to-solution}
\end{center}
\end{figure}

Figure \ref{fig:results/BIS/time-to-solution} shows the time-to-solution for the baroclinic instability test at $104 \times 1.25$ km resolution for set 2NC on 512 AMD Epyc (Rome) cores. All simulations use the 4th order tensor viscosity with non-dimensional parameters $(c_1,c_2,c_3)=(0.0045,0.0045,0.0001)$ with the viscosity computed at each stage of the IMEX method.  Each IMEX method uses the optimal configuration that admits the maximum usable time-step. We define the maximum usable time-step as the time-step that (i) yields the expected solution of the test case at 10 days and (ii) remains stable for 100 days. The figure shows that NHEVI-LU is 2x faster than NHEVI-GMRES and that LHEVI is much faster than either NHEVI method; 5x faster than NHEVI-LU and 10x faster than NHEVI-GMRES. In addition, the results show that ARK2 and ARS3 are always the fastest methods. 

\begin{table}
\begin{center}
\footnotesize
\begin{tabular}{ |l | l | r | r | r |}
\hline
Time-Integrator & HEVI Method &  $\Delta t$ (secs) & Update $\Delta t$ & Time-to-Solution (secs) \\ 
\hline \hline
ARK2 &  NHEVI-GMRES & 135&  -   & 1793 \\
ARS343 &  NHEVI-GMRES & 195 & -   & 1884 \\
ARK4 &  NHEVI-GMRES & 205 & -   & 2919 \\
ARK5 &  NHEVI-GMRES & 190 & -   & 4399\\
\hline
ARK2 &  NHEVI-LU         &  135 & -    & 502 \\
ARS343 &  NHEVI-LU         & 195 & -    & 515 \\
ARK4 &  NHEVI-LU         & 205 & -    & 791 \\
ARK5 &  NHEVI-LU         & 190 & -    & 1210 \\
\hline
ARK2 &  LHEVI                & 135 & 50 & 163 \\
ARS343 &  LHEVI               & 195 & 20    & 167 \\
ARK4 &  LHEVI               & 205 & 50    & 245 \\
 ARK5 &  LHEVI               & 190 & 50    & 373 \\
 \hline
\end{tabular}
\end{center}
 \caption{Baroclinic Instability. Configurations for the timing comparisons in Fig.\  \ref{fig:results/BIS/time-to-solution} for set 2NC. NHEVI uses $\epsilon=0.05$, with Newton tolerance $10^{-5}$ and GMRES tolerance $10^{-9}$.}
  \label{table:results/BIS/time-to-solution}
\end{table}

To help make sense of the time-to-solution for LHEVI, it is convenient to look at the ratio $R$ defined in \eqref{eq:time_integration/R_ratio}. Comparing ARK2 to the other methods we see that the ratio is $R_{2 \rightarrow 3}=1.04$, $R_{2 \rightarrow 4}=1.31$, $R_{2 \rightarrow 5}=2.48$, where the subscripts $2,3,4,5$ denote the order of the ARK methods; this ratio confirms that ARK2 and ARS3 result in a similar performance.  The cost of all the methods can be reduced if we only apply hyper-diffusion at the final stage. This reduces both the computational and communication costs substantially in a per time-step basis. We save this option for future work where we shall discuss the role of this strategy within multi-rate HEVI methods (see, e.g., \cite{mugg:2021}).

Now that we have compared the performance of each of the three HEVI methods and shown that LHEVI is, by far, the most efficient, we next need to quantify the accuracy of LHEVI compared to NHEVI. We could use this test case (baroclinic instability) for the time-convergence study but we choose a different case (inertia-gravity) presented below.

\subsection{Inertia-Gravity Wave}
\label{sec:inertia-gravity}

This section investigates the accuracy of the proposed LHEVI and NHEVI methods using the inertia-gravity wave from \cite{tomita:2004} (case 2-\#2). 
We use this test case in keeping with recent papers that use the two-dimensional inertia-gravity wave in a square domain \cite{skamarock:1994,baldauf:2016, baldauf:2021, waruszewski:2022}; in our case, we use its extension to three-dimensional global spherical domains suitable for verifying global atmospheric simulations.
 
\begin{figure}[h!]
	\centering
	\includegraphics[width=0.75\textwidth]{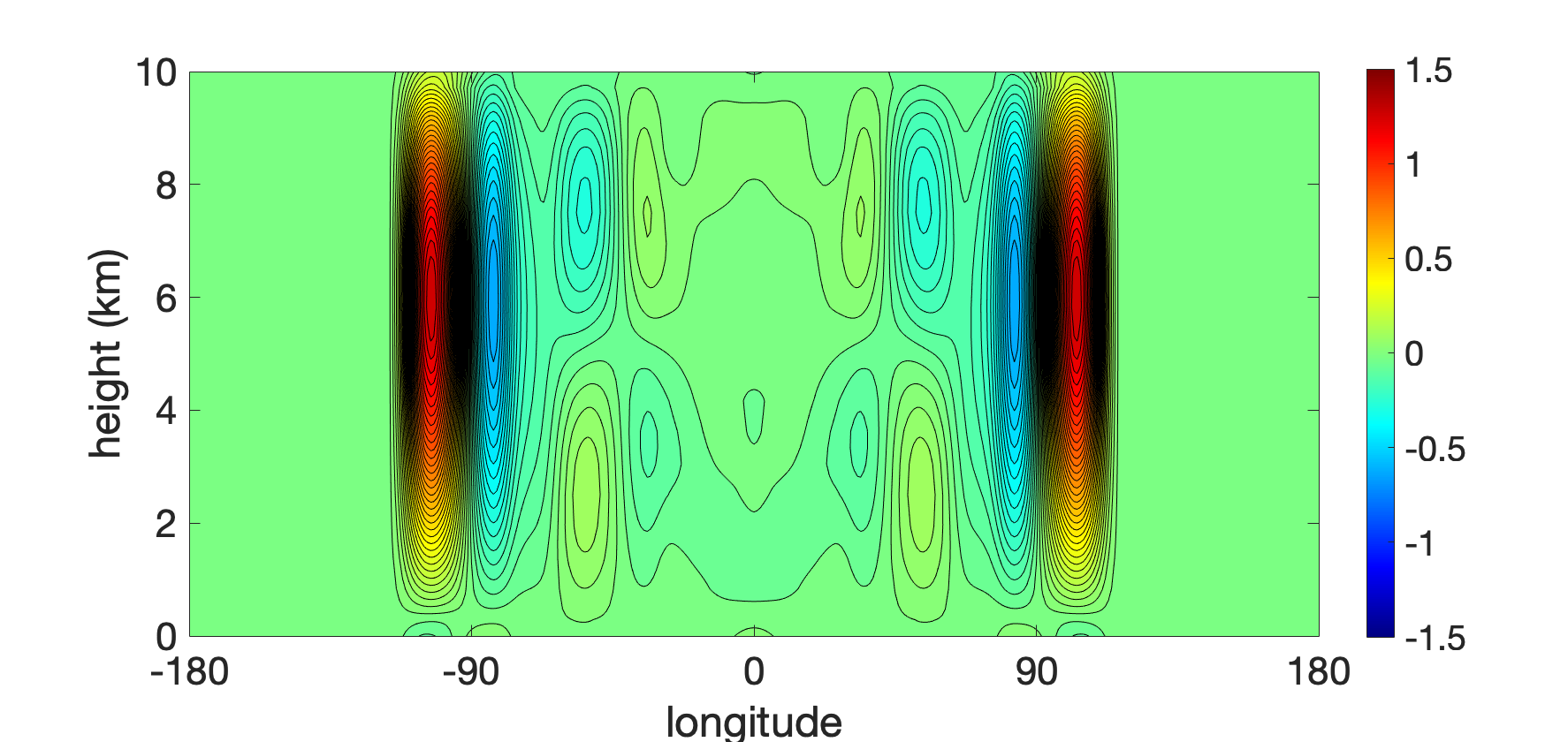}
	\caption{Inertia-gravity wave. Contours of potential temperature perturbation on the equator ($\phi=0)$ at $t=48$ hours using set 2NC with ARK5-LHEVI.}
  \label{fig:results/IGW/Theta}
\end{figure}

The wave propagation is initiated by adding the following potential temperature perturbation defined in spherical coordinates $(\lambda,\phi,z)$
\begin{equation}
  \theta' = \Delta\theta f(\lambda,\phi)g(z),
\end{equation}
where $\Delta\theta=10$ K is the amplitude of the potential temperature perturbation. The horizontal and vertical distribution functions, $f(\lambda,\phi)$ and $g(z)$, are defined as
\begin{align}
  f(\lambda,\phi) &= \left\{ 
    \begin{matrix}
      \frac{1}{2}\qty\big(1+\cos(\pi r/r_{pert})) & r<r_{pert},\\
      0 & r\ge r_{pert},
    \end{matrix}
  \right. \\
  g(z) &= \sin\qty(\dfrac{\pi z}{z_{top}})
\end{align}
where 
\begin{equation}
  r = r_e \cos^{-1}\qty\big(\sin\phi + \cos\phi_0 \cos\qty(\lambda))  
\end{equation}
with $z_{top}=10$ km, $r_{pert}=r_e/3$ km, and $\lambda$ and $\phi$ are the longitude and latitude, respectively.
The global domain is discretized using 4th-order polynomials in all three directions. The grid resolution in the horizontal and vertical directions are approximately 104 km and 417 m, respectively; this grid resolution corresponds to the same polynomial order and number of elements used in the baroclinic instability.  For temporal refinement, we test a series of time-step sizes, $\Delta t = 100\times 2^n$ seconds, with $n \in [-9,-8,\ldots,0]$. 
Figure \ref{fig:results/IGW/Theta} displays the potential temperature perturbation along the equator at the final time $t= $ 48 hr, calculated using ARK5-LHEVI.
The potential temperature perturbations qualitatively agree with Fig.\ 5b in \cite{tomita:2004}, indicating that the time-integrator can handle horizontally propagating gravity waves.

\begin{figure}[h!]
	\centering
  \begin{subfigure}{0.48\textwidth}
    \includegraphics[width=\textwidth]{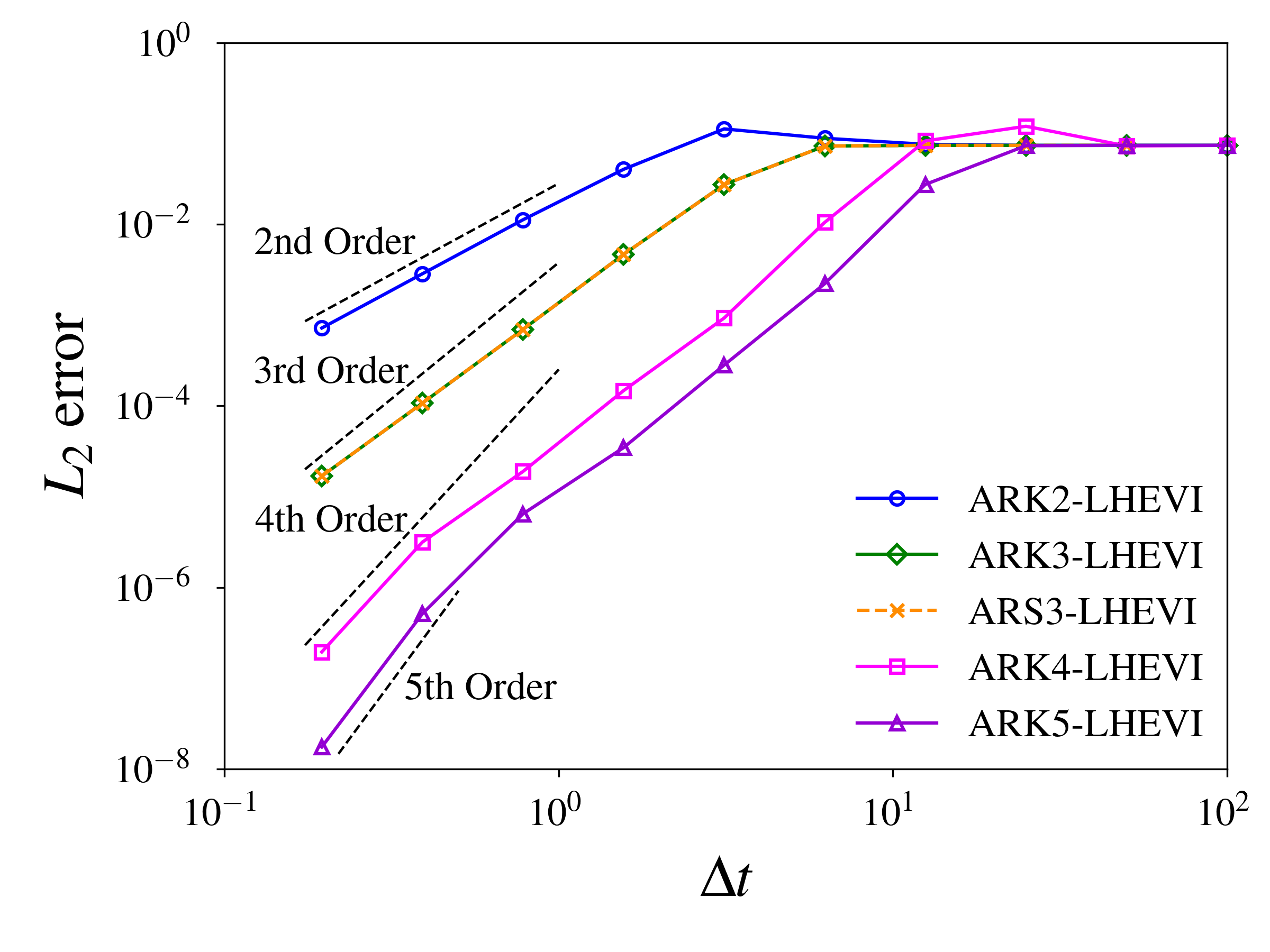}
    \caption{LHEVI}
  \end{subfigure}
  \begin{subfigure}{0.48\textwidth}
    \includegraphics[width=\textwidth]{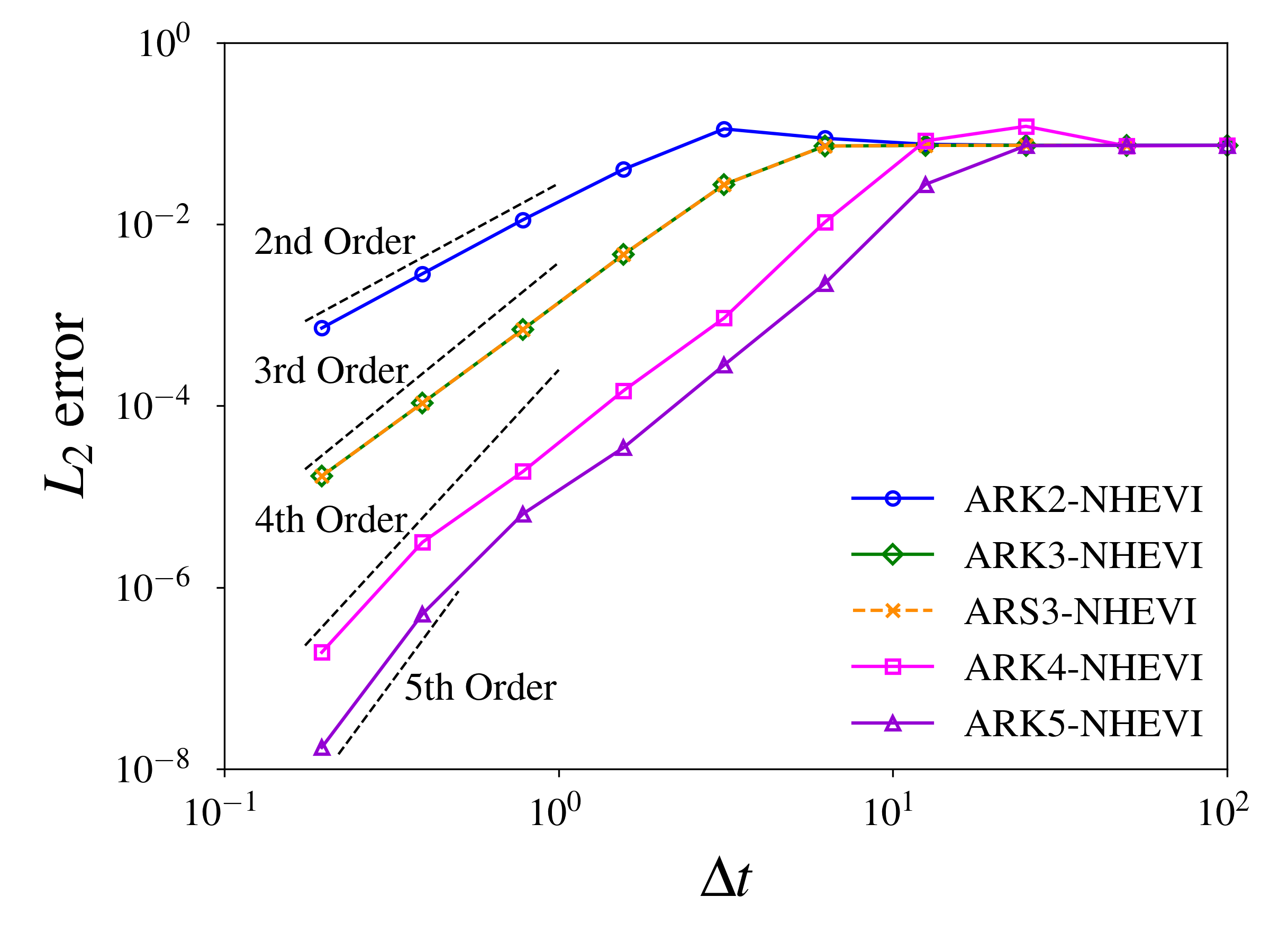}
    \caption{NHEVI}
  \end{subfigure}
	\caption{Inertia-gravity wave. Convergence of the normalized $L_2$-norm of the relative error with respect to time-step $\Delta t$ in the potential temperature field for LHEVI and NHEVI with set 2NC.}
  \label{fig:results/IGW/Convergence}
\end{figure}

We assess the accuracy of the method using a normalized $L_2$-norm of the relative error in the potential temperature field, which is defined as
\begin{equation}
  L_2\text{ error} = \dfrac{ \sqrt{\int_\Omega \left( \theta - \theta_{true} \right)^2 \, d\Omega}}{\sqrt{\int_\Omega \left( \theta_{true} \right)^2 \, d\Omega}}.
\end{equation}
In this test, we consider the solution obtained using a five-stage fourth-order low-storage RK method \cite{kennedy:2000} with $\Delta t = 100\times 2^{-12}$ seconds as the true solution and measure the relative errors at $t=3600$ seconds. 
Figure \ref{fig:results/IGW/Convergence} shows that the LHEVI and NHEVI methods achieve optimal convergence rates in the $L_2$-norm; the dashed lines denote the optimal convergence rates, which have slopes equal to their order on a log-log plot. Results are shown for set 2NC but similar results are obtained with the other two equation sets (not shown).  In the figure, the curves from ARK3 and ARS3 overlap, since both produce almost identical results. The optimal convergence rates demonstrate that the formulations of the LHEVI and NHEVI methods are consistent and yield accurate solutions; in fact, the accuracy of the LHEVI and NHEVI solutions are almost identical.
  The plateauing of the error norms for large $\Delta t$ was discussed in \cite{giraldo:2013} and is due to the acoustic waves being entirely damped.

In Sec.\ \ref{sec:baroclinic_instability}, we used the tensor-based hyper-diffusion to ensure the long-time stability of the methods. While the artificial hyper-diffusion parameters are designed to automatically scale with the grid and time-step sizes, as shown in \eqref{eq:viscosity_scaled}, for this particular test, we opt for the constant viscosities, $\nu_1=\nu_2=5.0\times 10^{7}$ and $\nu_3=150$, to eliminate the adaptive effects of time refinement and reach convergence (recall that to recover the physical viscosities requires squaring $\nu_i$).

\section{Conclusions}
\label{sec:conclusions}

We compared three versions of HEVI methods (NHEVI-GMRES, NHEVI-LU, and LHEVI) using a number of different IMEX methods found in the literature ranging from second to fifth order in accuracy. We confirmed that all methods obtained the expected rate of convergence for  LHEVI and NHEVI-LU. The LHEVI method performed the best regarding time-to-solution while achieving similar accuracy and stability characteristics to the NHEVI methods; we confirm these results by timing the simulations but also confirm the numerical results with the complexity analysis for each of the three HEVI methods. We tested the methods on the baroclinic instability on the sphere using 100-day simulations to confirm that the configurations for each HEVI method and all time-integrators were stable for long-time simulations. 
The baroclinic instability and inertia-gravity wave simulations indicate that the optimal choice of time-integrator is LHEVI with either ARK2 or ARS3, as both methods yield similar time-to-solution and relative $L_2$ error at their maximum usable time-steps.  All time-integrators conserve mass with any choice of equation set; in addition, total energy is conserved up to dissipation due to hyper-diffusion for set 3C.
In the future, we will test LHEVI with ARK2 and ARS3 on more stringent flow problems that produce steeper vertical gradients to test whether the same conclusions hold. However, the results of this study offer new possibilities for constructing yet more efficient time-integrators, e.g., extending the HEVI methods to fully implicit methods in all dimensions that then employ the corresponding version of LHEVI as preconditioners. We will report on this class of methods in the future and test them on more realistic NWP and space weather problems.

\appendix

\section{Sets 2C and 3C}
\label{sec:appendix/set2c_set3c}

For completeness, below we include the governing equations for sets 2C and 3C in \ref{sec:appendix/set2c_set3c/governing_equations}, the horizontal and vertical operators in HEVI in \ref{sec:appendix/set2c_set3c/hevi_operators}, and the HEVI operator Jacobians in \ref{sec:appendix/set2c_set3c/operator_jacobian}.

\subsection{Governing Equations}
\label{sec:appendix/set2c_set3c/governing_equations}

A simple manipulation of \eqref{eq:set2nc} results in the conservation form
\begin{subequations}
\label{eq:set2c}
\begin{equation}
\label{eq:set2c_mass}
\diff{\rho}{t} + \nabla \cdot \vc{U} = 0
\end{equation}
\begin{equation}
\label{eq:set2c_momentum}
\diff{\vc{U}}{t} + \nabla \cdot \left( \frac{\vc{U}\otimes \vc{U}}{\rho} + P \vc{I}_d \right)  + \rho \nabla \Phi +  2 \omega \hat{\vc{r}} \times \vc{U} = \vc{0}
\end{equation}
\begin{equation}
\label{eq:set2c_energy}
\diff{\Theta}{t} + \nabla \cdot \left( \frac{\Theta \vc{U}}{\rho} \right)  = 0
\end{equation}
\begin{equation}
\label{eq:set2c_pressure}
P = P_A \left( \frac{R \Theta}{P_A} \right)^{\gamma},
\end{equation}
\end{subequations}
where $\vc{U}=\rho \vc{u}$ is momentum, $\Theta=\rho \theta$ is density potential temperature, and $\vc{I}_d$ is the rank-$d$ identity matrix where $d$ is the spatial dimension; with the proper numerics (such as the EBG methods we use in this work) mass and integrated potential temperature are globally conserved.  Set 2C, or a variant thereof, is used in, e.g., ICON \cite{zangl:2015} and MPAS \cite{skamarock:2012}.

Set 3C is the most commonly used form of the compressible Euler equations in computational fluid dynamics (CFD) and is written as follows
\begin{subequations}
\label{eq:set3c}
\begin{equation}
\label{eq:set3c_mass}
\diff{\rho}{t} + \nabla \cdot \vc{U} = 0
\end{equation}
\begin{equation}
\label{eq:set3c_momentum}
\diff{\vc{U}}{t} + \nabla \cdot \left( \frac{\vc{U}\otimes \vc{U}}{\rho} + P \vc{I}_d \right)  + \rho \nabla \Phi +  2 \omega \hat{\vc{r}} \times \vc{U} = \vc{0}
\end{equation}
\begin{equation}
\label{eq:set3c_energy}
\diff{E}{t} + \nabla \cdot \left( \frac{(E+P) \vc{U}}{\rho} \right)  = 0
\end{equation}
\begin{equation}
\label{eq:set3c_pressure}
P = \gmo \left( E - \frac{\vc{U} \cdot \vc{U} }{2 \rho} - \rho \Phi \right),
\end{equation}
\end{subequations}
where $E=\rho e$ is total energy, with $e=c_v T + \frac{\vc{u} \cdot \vc{u}}{2} + \Phi$ denoting the total specific energy (internal, kinetic, and potential).  This set will conserve both mass and total energy. Unlike sets 2NC and 2C, 3C does not make any assumptions on the composition of the atmosphere and is valid for a whole atmosphere model.
Set 3C, to our knowledge, has not been used in any atmospheric model except for NUMA \cite{giraldo:2008a,giraldo:2010}, an early version of CLIMA \cite{sridhar:2022}, and Atum \cite{waruszewski:2022,souza:2023}. 

\subsection{HEVI Operators}
\label{sec:appendix/set2c_set3c/hevi_operators}

The horizontal and vertical HEVI operators for set 2C are
\begin{equation}
  \vc{H}(\vc{q}) =  - \begin{pmatrix}
       \frac{1}{J} \nabla_{H} \cdot \left( J \vc{F}_{\rho}^{H} \right) \\
        \frac{1}{J} \nabla_{H} \cdot \left( J \vc{F}_{{U}}^{H} \right) +  \rho \left(\diff{\Phi}{\xi}\nabla \xi + \diff{\Phi}{\eta}\nabla \eta \right) + 2\omega \wh{\vc{\zeta}} \times \vc{U}  \\
     \frac{1}{J} \nabla_{H} \cdot \left( J \vc{F}_{T}^{H} \right) \\
  \end{pmatrix} + \mathcal{H}_{\nu}(\qvector)
\end{equation}
\begin{equation}
  \vc{V}\left(\vc{q}\right) = - \begin{pmatrix}
    \frac{1}{J} \diff{ }{\zeta} \left( J \vc{F}_{\rho}^\zeta \right) \\
    \frac{1}{J} \diff{ }{\zeta} \left( J \vc{F}_{U}^\zeta \right) + \rho \diff{\Phi}{\zeta}\nabla\zeta \\
     \frac{1}{J} \diff{ }{\zeta} \left( J \vc{F}_{T}^\zeta \right) 
  \end{pmatrix}
\end{equation}
where $\vc{F}_{T}$ denotes the flux for the thermodynamic variable and the contravariant flux components are defined as
\( \vc{F}_q^{\xi,\eta} = \left(\vc{F}_q^\xi, \vc{F}_q^\eta, 0\right)^\mathcal{T} \), $\vc{F}_q^{\zeta}$,  for the state vector $\qvector=(\rho, \vc{U}\transpose, \Theta)\transpose$, with the contravariant fluxes computed as \( \vc{F}_q^{i} = \vc{F}_q \cdot \nabla \xi^i \).
Set 3C has a similar form except that the thermodynamic variable and corresponding flux are different, along with the equation of state (compare Eqs.\ \eqref{eq:set2c} and \eqref{eq:set3c}).

\subsection{HEVI Operator Jacobians}
\label{sec:appendix/set2c_set3c/operator_jacobian}
 
For set 2C we arrive at the following operator Jacobian
\[
\vc{\mathbb{J}}_{2C}(\vc{q})=\vc{\mathbb{I}} + \Lambda \mathbb{D}^{\zeta} \vc{\mathbb{K}}_{2C}(\vc{q}_j) +  \Lambda \vc{\mathbb{P}}
\]
where $\mathbb{D}^{\zeta}$ can be either $\mathbb{D}^{\zeta}_{strong}$ or $\mathbb{D}^{\zeta}_{weak}$, denoting the strong or weak forms, respectively. These are defined as follows
\begin{subequations}
\label{eq:appendix/differentiation_matrix}
\be
\label{eq:appendix/differentiation_matrix/strong}
\left( \mathbb{D}^{\zeta}_{strong} \right)_{ij} = \frac{J_j}{J_i} D^{\zeta}_{ij}
\ee
\be
\label{eq:appendix/differentiation_matrix/weak}
\left( \mathbb{D}^{\zeta}_{weak} \right)_{ij} = -  \frac{(\varpi J)_j}{(\varpi J)_i} D^{\zeta}_{ji},
\ee
\end{subequations}
where
{\scriptsize
\begin{equation}
\label{eq:set2c_jacobian}
 \vc{\mathbb{K}}_{2C}(\vc{q}_j) = 
 \begin{pmatrix}
  0 & 
 \zeta_x &
 \zeta_y &
\zeta_z &
 0 \\
 -\frac{U}{\rho^2} U^{\zeta} + \zeta_x \diff{P}{\rho} & 
 \frac{1}{\rho} U^{\zeta} + \zeta_x \left( \frac{U}{\rho} + \diff{P}{U} \right) &
 \frac{1}{\rho} \zeta_y \frac{U}{\rho} + \zeta_x \diff{P}{V} &
\frac{1}{\rho} \zeta_z \frac{U}{\rho} + \zeta_x \diff{P}{W}  &
 \frac{1}{\rho} \zeta_x \diff{P}{E} \\ 
 -\frac{V}{\rho^2} U^{\zeta} + \zeta_y \diff{P}{\rho} & 
  \frac{1}{\rho} \zeta_x \frac{V}{\rho} + \zeta_y \diff{P}{U} &
  \frac{1}{\rho} U^{\zeta} + \zeta_y \left( \frac{V}{\rho} + \diff{P}{V} \right) &
\frac{1}{\rho} \zeta_z \frac{V}{\rho} + \zeta_y \diff{P}{W}  &
 \frac{1}{\rho} \zeta_y \diff{P}{E} \\ 
 -\frac{W}{\rho^2} U^{\zeta} + \zeta_z \diff{P}{\rho} & 
  \frac{1}{\rho} \zeta_x \frac{W}{\rho} + \zeta_z \diff{P}{U} &
 \frac{1}{\rho} \zeta_y \frac{W}{\rho} + \zeta_z \diff{P}{V}  &
  \frac{1}{\rho} U^{\zeta} + \zeta_z \left( \frac{W}{\rho} + \diff{P}{W} \right) &
 \frac{1}{\rho} \zeta_z \diff{P}{E} \\ 
- \frac{U^{\zeta}}{\rho^2} & 
\zeta_x \frac{\Theta}{\rho}  &
\zeta_y \frac{\Theta}{\rho}  &
\zeta_z \frac{\Theta}{\rho}  &
\frac{U^{\zeta}}{\rho} & 
 \end{pmatrix}
\end{equation}
}
\noindent with all terms evaluated at $j$, and 
\begin{equation}
\label{eq:set2c_jacobian_part2}
 \vc{\mathbb{P}} = 
 \mbox{diag} \begin{pmatrix}
  0 \\ 
 \zeta_x \diff{\Phi}{\zeta} \\
\zeta_y \diff{\Phi}{\zeta} \\
 \zeta_z \diff{\Phi}{\zeta}  \\
 0
\end{pmatrix},
\end{equation}
with all terms evaluated at $i$.

For set 3C we get
\[
\mathbb{J}_{3C}(\vc{q})=\vc{\mathbb{I}} + \Lambda \mathbb{D}^{\zeta} \mathbb{K}_{3C}(\vc{q_j})  +  \Lambda \vc{\mathbb{P}}
\]
with $\mathbb{D}^{\zeta}$ defined as in Eq.\ \eqref{eq:appendix/differentiation_matrix} and $\vc{\mathbb{P}}$ as in Eq.\ \eqref{eq:set2c_jacobian_part2}, with
{\scriptsize
\begin{equation}
\label{eq:set3c_jacobian}
 \mathbb{K}_{3C}(\vc{q_j}) = 
 \begin{pmatrix}
  0 & 
 \zeta_x  &
 \zeta_y  &
\zeta_z  &
 0 \\
 -\frac{U}{\rho^2} U^{\zeta} + \zeta_x \diff{P}{\rho} & 
 \frac{1}{\rho} U^{\zeta} + \zeta_x \left( \frac{U}{\rho} + \diff{P}{U} \right) &
 \frac{1}{\rho} \zeta_y \frac{U}{\rho} + \zeta_x \diff{P}{V} &
\frac{1}{\rho} \zeta_z \frac{U}{\rho} + \zeta_x \diff{P}{W}  &
 \frac{1}{\rho} \zeta_x \diff{P}{E} \\ 
 -\frac{V}{\rho^2} U^{\zeta} + \zeta_y \diff{P}{\rho} & 
  \frac{1}{\rho} \zeta_x \frac{V}{\rho} + \zeta_y \diff{P}{U} &
  \frac{1}{\rho} U^{\zeta} + \zeta_y \left( \frac{V}{\rho} + \diff{P}{V} \right) &
\frac{1}{\rho} \zeta_z \frac{V}{\rho} + \zeta_y \diff{P}{W}  &
 \frac{1}{\rho} \zeta_y \diff{P}{E} \\ 
 -\frac{W}{\rho^2} U^{\zeta} + \zeta_z \diff{P}{\rho} & 
  \frac{1}{\rho} \zeta_x \frac{W}{\rho} + \zeta_z \diff{P}{U} &
 \frac{1}{\rho} \zeta_y \frac{W}{\rho} + \zeta_z \diff{P}{V}  &
  \frac{1}{\rho} U^{\zeta} + \zeta_z \left( \frac{W}{\rho} + \diff{P}{W} \right) &
 \frac{1}{\rho} \zeta_z \diff{P}{E} \\ 
U^{\zeta} \left( \frac{1}{\rho} \diff{P}{\rho} - \frac{E+P}{\rho^2} \right) & 
\frac{U^{\zeta}}{\rho} \diff{P}{U} + \frac{E+P}{\rho} \zeta_x &
 \frac{U^{\zeta}}{\rho} \diff{P}{V} + \frac{E+P}{\rho} \zeta_y &
  \frac{U^{\zeta}}{\rho} \diff{P}{W} + \frac{E+P}{\rho} \zeta_z &
  \frac{U^{\zeta}}{\rho} \left( 1 + \diff{P}{E} \right)
  \end{pmatrix}.
\end{equation}
}
with all terms evaluated at $j$.

Upon differentiation by $\mathbb{D}^{\zeta} $, numerical integration (multiplying by $\varpi_i J_i$), and applying DSS along $\zeta$, results in the global Jacobian matrices $\mathcal{G} \left( \vc{\mathbb{J}} \right)
\in \mathbb{R}^{M \times M}$, where $M=n_{var}N^{\zeta}_e\left( N_{\zeta}+1 \right)$, and $\mathcal{G}$ represents the operator that takes ${\mathbb{J}}$ and constructs its global representation, as presented in Sec.\ \ref{sec:nhevi-lu}.

\section*{Acknowledgments}  \label{sct:acknowledgement}
Felipe Alves was funded by the Defense Sciences Office of the Defense Applied Research and Projects Agency (DARPA
DSO) through the Space Environment Exploitation (SEE) program.  F.\ X.\ Giraldo, Soonpil Kang, James F.\ Kelly, and P.\ Alex Reinecke gratefully acknowledge the support of the Office of Naval Research under grant \# N0001419WX00721. F.\ X.\ Giraldo was also supported by the National Science Foundation under grant AGS-1835881. This work was performed when Felipe Alves and Soonpil Kang held National Academy of Sciences’ National Research Council (NRC) Fellowships at the Naval Postgraduate School.

\bibliographystyle{siam}
\bibliography{bibliography/Giraldo_refs_v3.bib,bibliography/NPS.bib}

\end{document}